\documentclass[12pt]{article}
\usepackage{paperstyle}
\usepackage{macros}

\begin{document}

\title{Stable closed geodesics and stable figure eights in convex hypersurfaces}
\author{Herng Yi Cheng}
\date{}

\maketitle

\begin{abstract}
    For each odd $n \geq 3$, we construct a closed convex hypersurface of $\mathbb{R}^{n+1}$ that contains a non-degenerate closed geodesic with Morse index zero. A classical theorem of J.\,L.~Synge would forbid such constructions for even $n$, so in a sense we prove that Synge's theorem is ``sharp.'' We also construct stable figure eights: that is, for each $n \geq 3$ we embed the figure eight graph in a closed convex hypersurface of $\mathbb{R}^{n+1}$, such that sufficiently small variations of the embedding either preserve its image or must increase its length.
    
    These stable geodesics and stable figure eights are mainly obtained by constructing explicit billiard trajectories with ``controlled parallel transport'' in convex polytopes.
\end{abstract}

2020 \emph{Mathematics Subject Classification.} Primary 53C22; Secondary 53C42

\section{Introduction}

\subsection{Stable closed geodesics}

Is it possible to lasso a convex body? Or, more formally, can a convex hypersurface of Euclidean space contain a \emph{stable closed geodesic}---a closed geodesic with Morse index zero? With some additional hypotheses, the answer will be ``no'': let $M$ be a Riemannian manifold with positive sectional curvature. If $M$ is even-dimensional and orientable, then a classical theorem of J.\,L.~Synge prevents it from having an stable closed geodesic \cite{Synge_GeodesicPosCurv}. If $M$ is simply connected and compact, and has curvature that is more than $\frac14$-pinched, then it also cannot contain a stable closed geodesic: W.\,P.\,A.~Klingenberg's injectivity radius estimate and the Morse-Schoenberg comparison theorem would imply that every closed geodesic in $M$ has index at least $\dim M - 1$ \cite{Klingenberg_RiemGeom}.

On the other hand, we prove that Synge's theorem is ``sharp'' in some sense by constructing stable closed geodesics in odd-dimensional convex hypersurfaces.

\begin{theorem}
    \label{thm:StableClosedGeodesicPosCurv}
    For every odd integer $n \geq 3$, there exists a closed convex hypersurface $M_n'$ of $\R^{n+1}$ with positive sectional curvature that contains a simple and stable closed geodesic.
    
    Furthermore, for each odd $n \geq 3$, every Riemannian $n$-sphere whose metric is sufficiently close to that of $M_n'$ in the $C^\infty$ topology also contains a simple and stable closed geodesic.
\end{theorem}

(\emph{Simple} means that the closed geodesic does not self-intersect.)

Our convex hypersurfaces $M_n'$ are simply connected and positively curved, which are in some sense the most difficult class of spaces to find stable closed geodesic in. For instance, in any closed Riemannian manifold that is not simply connected, the shortest non-contractible loop is a stable closed geodesic. In addition, any closed geodesic whose neighbourhood is negatively curved has to be stable, by the formula for the second variation of length.

If we consider closed geodesics as analogues of minimal surfaces, then our result contrasts with the fact that closed Riemannian manifolds with positive Ricci curvature cannot contain stable, two-sided, and closed minimal hypersurfaces \cite[p.~79]{Simons_MinimalVarieties}.

We construct each stable closed geodesic in an ``elementary'' way using line segments in the smooth part of a manifold that has a flat Riemannian metric away from some singularities; we will eventually smooth these manifolds into convex hypersurfaces. Indeed, the geodesics involved are derived from certain billiard trajectories in convex polytopes that can be explicitly computed.

It should be highlighted that W.~Ziller constructed a 3-dimensional positively-curved homogeneous space that is diffeomorphic to $S^3$ and contains a stable closed geodesic \cite[Example~1]{Ziller_HomogSpaceStableClosedGeodesic}. It may be possible that Ziller's construction can be generalized to find stable closed geodesics in positively-curved homogeneous spaces of higher odd dimension. However, these spaces could never be isometric to closed convex hypersurfaces of Euclidean space; if that were possible, then by a result of S.~Kobayashi \cite{Kobayashi_HomogHypersurfaces} they would have to be isometric to round spheres, which cannot have stable closed geodesics.

\subsection{Stable figure eights}

Our other main result is the construction of \emph{stable figure eights} in convex hypersurfaces. Stable figure eights are roughly embeddings of the figure eight graph into a Riemannian manifold so that small perturbations of the embedding must increase its length. To define it, let us first introduce \emph{stationary geodesic nets} in a Riemannian manifold $M$, which are immersions in $M$ of a connected graph $\Gamma$ that are stationary points in the space of such immersions under the length functional. That is, every perturbation of the immersion has a vanishing first variation in length.  Geometrically speaking, the edges of stationary geodesic nets are geodesics, and each vertex satisfies the \emph{stationarity condition}: the unit tangent vectors of the adjacent edges that point away from the vertex must sum to zero. Stationary geodesic nets can also be defined as stationary 1-dimensional integral varifolds in geometric measure theory \cite{Pitts_StationaryVarifolds,AllardAlmgren_StationaryVarifolds}. Thus stationary geodesic nets can be regarded as generalizations of closed geodesics and 1-dimensional analogues of minimal surfaces, and studying their properties, such as their stability and Morse indices, could inform our understanding of closed geodesics and minimal surfaces.

There have been several results proving the existence of stationary geodesic nets in positively-curved Riemannian 2-spheres by H.~Morgan \cite{HassMorgan_GeodesicNetS2}, and in general closed Riemannian manifolds by A.~Nabutovsky and R.~Rotman \cite{NabutovskyRotman_Stationary1Cycle,NabutovskyRotman_GeodesicNets}, and Rotman \cite{Rotman_ShortestGeodesicNet,Rotman_Flowers}. Y.~Liokumovich and B.~Staffa even proved that the union of stationary geodesic nets is dense in a generic closed Riemannian manifold  \cite{LiokumovichStaffa_GenericDensity}. On the other hand, none of these results are able to guarantee that the resulting stationary geodesic nets contain closed geodesics or do not contain them, with the exception of a result of Morgan that Riemannian 2-spheres sufficiently close to the round sphere contain stationary geodesic nets shaped like $\theta$ graphs (two vertices connected by three edges) \cite[p.~3848]{HassMorgan_GeodesicNetS2}. It would be interesting to prove that manifolds from some broad class contain stationary geodesic nets that are \emph{irreducible}, in the sense that no pair of tangent vectors at the basepoint are parallel, as that would exclude ``degenerate'' stable geodesic bouquets that are merely unions of stable closed geodesics. This would help establish limits on how much techniques that produce stationary geodesic nets (e.g. Almgren-Pitts min-max theory, see the survey \cite{Marques_MinimalSurfacesSurvey}) can also help to produce closed geodesics.

A \emph{stable geodesic net} is a stationary geodesic net that must lengthen under any sufficiently small perturbation. (It will be formally defined in \cref{sec:Definitions}.) We will mostly consider stable geodesic nets that are immersions of ``flower-shaped'' graphs $\bq_k$ comprising $k \geq 1$ loops at a single vertex. Such nets are called \emph{stable geodesic bouquets}. When $k = 2$, such nets are called \emph{stable figure eights}. Figure eights are the simplest shape of geodesic net, hence stable figure eights may be viewed as the simplest generalizations of stable closed geodesics. 

Like stable closed geodesics, the existence of stable geodesic nets can be obstructed by certain positive curvature conditions. A well-known conjecture of H.\,B.~Lawson and J.~Simons \cite{LawsonSimons} would imply that if $M$ is closed and simply connected, and its curvature is $\frac14$-pinched, then it cannot contain any stable submanifolds or stable varifolds of any dimension, which includes stable closed geodesics and stable geodesic bouquets. This conjecture has been proven under various additional pinching assumptions \cite{Howard_SuffPinched,ShenXu_QuarterPinchedHypersurfaces,HuWei_FifthPinched}. In particular, the conjecture holds for all round spheres, so that none of them contain any stable geodesic net. 

I.~Adelstein and F.\,V.~Pallete also proved that positively curved Riemannian 2-spheres cannot have a stable geodesic net in the shape of a $\theta$ graph \cite{AdelsteinPallete}. This allowed them to prove that the length of the shortest closed geodesic in a positively curved Riemannian 2-sphere $M$ is at most $3d$, where $d$ is the diameter of $M$. This sharpened the upper bound of $4d$ proved by A.~Nabutovsky and R.~Rotman \cite{NabutovskyRotman_ShortestClosedGeodesicS2} and independently by S.~Sabourau \cite{Sabourau_ShortestClosedGeodesic} for general Riemannian 2-sphere. This improvement was possible by showing that certain processes of contracting graphs embedded in $M$ cannot get stuck in stable $\theta$ graphs, and will converge to stable closed geodesics. This suggests that the existence of stable geodesic nets may have implications on estimates of the lengths of closed geodesics.

In negatively curved manifolds, stable geodesic nets are easier to construct, akin to stable closed geodesics. If a stationary geodesic net lies in a region of strictly negative curvature then it must also be stable due to the formula for the second variation of length.

Given the negative results in positive curvature, it is natural to ask which, if any, Riemannian spheres with positive sectional curvature can contain stable geodesic bouquets. The author recently constructed, in every dimension $n \geq 2$, a Riemannian $n$-sphere $S_n$ which is isometric to a convex hypersurface and which contains an irreducible stable geodesic bouquet. The bouquet has 3 loops if $n = 2$, and has $n$ loops if $n \geq 3$ \cite{Cheng_StableGeodesicNets}. The convex hypersurfaces have strictly positive sectional curvature. To our understanding, that was the first existence result for irreducible stable geodesic nets in a simply connected and closed manifold with strictly positive sectional curvature.

In this paper we prove another existence result in which only two loops are necessary in the stable geodesic bouquet. Irreducible stable figure eights are the simplest possible stable geodesic nets that contain no closed geodesics.

\begin{theorem}
    \label{thm:Stable2LoopPosCurv}
    For every integer $n \geq 3$, there exists a closed convex hypersurface $M_n$ of $\R^{n+1}$ with positive sectional curvature that contains a simple, irreducible, and stable figure eight.
    
    Furthermore, for each $n \geq 3$, every Riemannian $n$-sphere whose metric is sufficiently close to that of $M_n$ in the $C^\infty$ topology also contains a simple, irreducible, and stable figure eight.
\end{theorem}

(By \emph{simple} we mean that the stable figure eight is an injective immersion. Note that every stable figure eight in a Riemannian surface (i.e. $n = 2$) must be the image of a self-intersecting closed geodesic, due to the stationarity condition.)

\Cref{thm:Stable2LoopPosCurv} shows that Synge's theorem is ``sharp'' in another sense: it becomes false when stable closed geodesics are replaced with their simplest generalizations among the geodesic nets, stable figure eights.

On the face of it, this result seems stronger than the one in \cite{Cheng_StableGeodesicNets} in the sense that when we are allowed to have more loops, it seems easier to ``arrange'' the loops into a stable geodesic bouquet. Using only two geodesic loops $\gamma_1, \gamma_2 : I \to M$ gives us much less freedom, because the condition of being a critical point of the length functional requires that the same line must bisect the angle between $\gamma_1'(0)$ and $-\gamma_1'(1)$ as well as the angle between $\gamma_2'(0)$ and $-\gamma_2'(1)$. Moreover, both of those angles must be equal in magnitude.

Nevertheless, \cref{thm:Stable2LoopPosCurv} is independent from the result in \cite{Cheng_StableGeodesicNets}. Using only two loops forces us to control the parallel transport map along the geodesic loops more precisely, and we achieve this using the novel techniques used to prove \cref{thm:StableClosedGeodesicPosCurv}. Our stable figure eights will, like our stable closed geodesics, be constructed in an elementary way from billiard trajectories in convex polytopes.

(A billiard trajectory in a convex polytope $\X \subset \R^n$ is a path that travels in a straight line until it collides with the interior of a face of $\X$, after which its velocity vector gets reflected about that face and the trajectory proceeds with the new velocity in another straight line, and so on.\footnote{For a formal definition, see \cite[Section~4.1]{Cheng_StableGeodesicNets}.})

\subsection*{Acknowledgements}

The author would like to express their gratitude to their academic advisors Alexander Nabutovsky and Regina Rotman for suggesting this research topic, and for valuable discussions. The author would also like to thank Isabel Beach for useful discussions. The author was supported by the Vanier Canada Graduate Scholarship.

\section{Summary of Key Ideas}
\label{sec:KeyIdeas}

\subsection{Stable closed geodesics as core curves of twisted tubes}
\label{sec:Intro_CoreCurves}

Synge proved that an even-dimensional and orientable manifold $M$ cannot contain a stable closed geodesic by demonstrating that any closed geodesic $\gamma$ in $M$ can be shortened by perturbing it along some parallel vector field orthogonal to $\gamma$. Thus our strategy is to construct odd-dimensional manifolds that contain a closed geodesic without parallel vector fields orthogonal to it. For instance, given some integer $n \geq 2$ and lengths $r, \ell > 0$, consider the space $(\D^{n-1}_r \times [0,\ell])/{\sim}$, where $\D^{n-1}_r$ is the closed disk of radius $r$ centered at the origin in $\R^{n-1}$, and the equivalence relation $\sim$ identifies $(x,0)$ with $(-x,\ell)$ for all $x \in \D^{n-1}_r$. The resulting space $T$ is a flat Riemannian manifold with boundary, which we call an \emph{$n$-dimensional twisted tube}. (When $n = 2$ this is a M\"obius strip.) $T$ is orientable if and only if $n$ is odd. Parallel transport along the core curve $\gamma$ negates every vector orthogonal to $\gamma$, so $\gamma$ has no parallel vector fields orthogonal to it. An elementary argument shows that $\gamma$ attains the minimal length among all smooth closed curves in $T$ freely homotopic to $\gamma$.\footnote{The curve with minimal length must be a geodesic. Otherwise the curve-shortening flow would reduce its length. (The convexity of $T$ implies that curves in $T$ will stay in $T$ under that flow.) However, the only geodesics that stay in $T$ are curves that run parallel to $\gamma$. Being homotopic to $\gamma$ means that the minimizing length curve must close up after exactly one round around $T$, but that can only happen if the curve passes through the origin of some copy of $\D^{n-1}_r$ in $T$. Thus that curve must be $\gamma$.} We will prove later that $\gamma$ is stable (see \cref{Lem:MaxTwistStableGeodesic}).

Twisted tubes will help us to prove \cref{thm:StableClosedGeodesicPosCurv} in the following way. For each odd $n \geq 3$ we will first embed an $n$-dimensional twisted tube in a ``polyhedral manifold'' obtained by gluing two copies of a convex $n$-dimensional polytope $\X^n \subset \R^n$ (henceforth, \emph{$n$-polytope}) along their boundaries via the identity map. The result, called the \emph{double} of $\X^n$ and denoted by $\dbl\X^n$, is a topological manifold homeomorphic to the $n$-sphere. In fact, $\dbl\X^n$ is a Riemannian manifold with singularities at the image of the $(n-2)$-skeleton of $\X^n$. The \emph{smooth part} of $\dbl\X^n$ is the complement of its singularities, denoted by $\dbl[sm]\X^n$, which has a flat Riemannian metric induced from the Euclidean metric on $\X^n$. The twisted tube will be isometrically embedded in $\dbl[sm]\X^n$, and its core curve will be a stable closed geodesic $\gamma_n$. The final step will then be to smooth $\dbl\X^n$ into a convex hypersurface in $\R^{n+1}$ while preserving the stability of $\gamma_n$, using arguments from \cite{Cheng_StableGeodesicNets}. Here it is important that $\gamma_n$ is non-degenerate, to ensure that after sufficiently small changes in the ambient metric from the smoothing, $\gamma_n$ will be close to a new stable closed geodesic.

\begin{figure}[p]
    \centering
    \begin{tabular}{c}
        \includegraphics[width=0.8\linewidth]{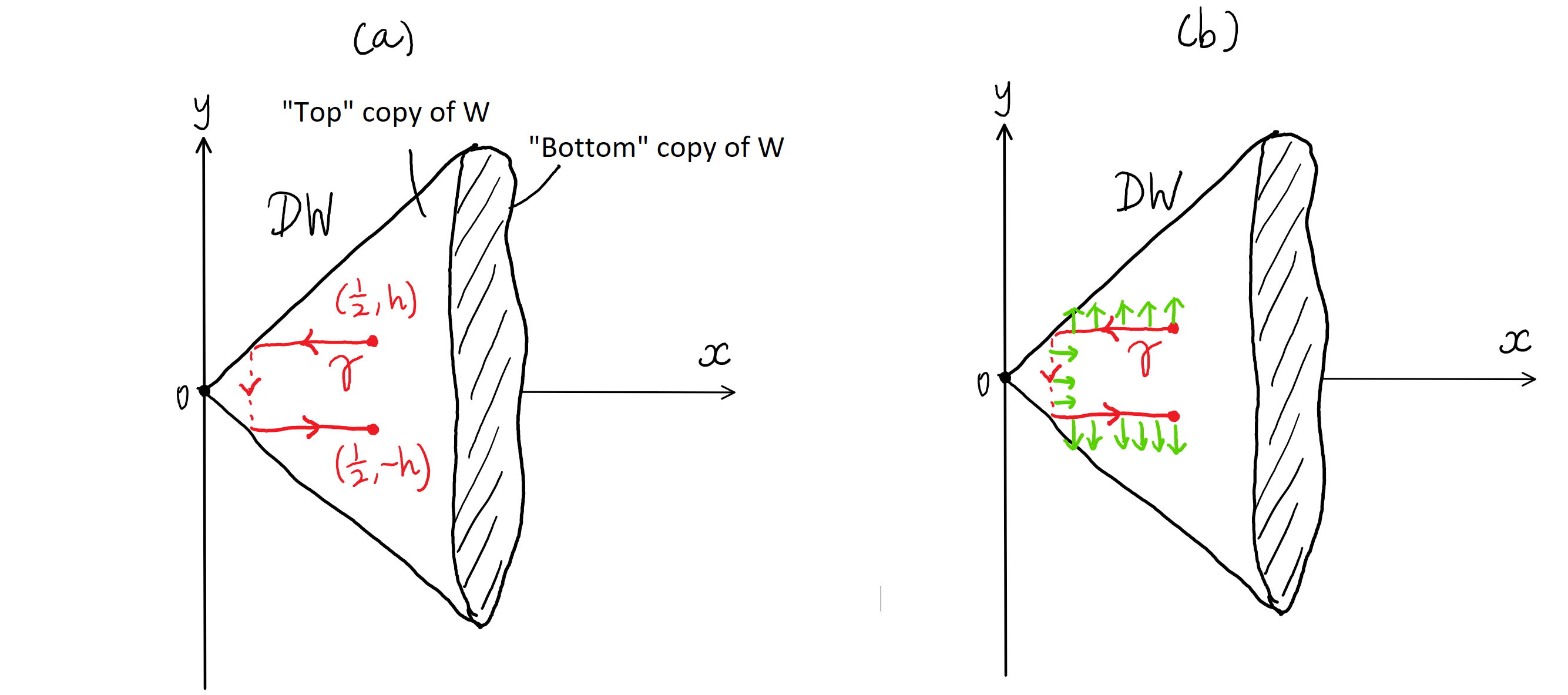}
    \\
        \includegraphics[width=0.8\linewidth]{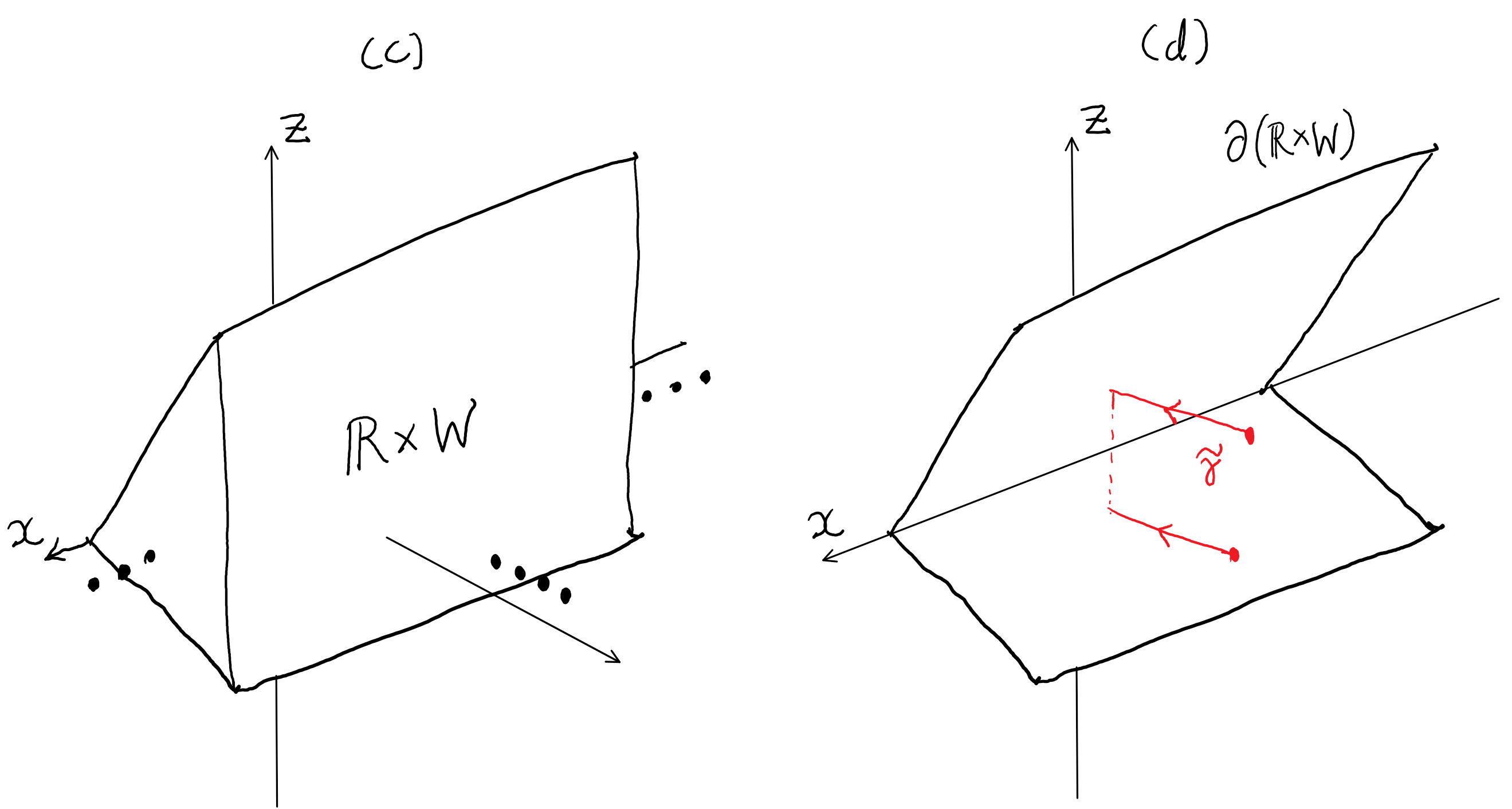}
    \end{tabular}
    \caption{(a) $\dbl\W$ is a ``cone,'' with both copies of $\W$ plotted on the $xy$-plane. $\gamma$ winds around the vertex of the cone, traveling from the ``top copy'' of $\W$ to the ``bottom copy'' and then back. (b) An illustration of the parallel transport of the green vectors. (c) An illustration of part of $\R \times \W$, which extends to infinity in three directions as indicated by the ellipses. (d) $\tilde\gamma$ travels in $\dbl(\R \times \W)$ through the depicted copy of $\R \times \W$, through the other copy (shown as a dashed segment), and then back to the first copy.}
    \label{fig:Wedges}
\end{figure}

The key idea behind the construction of $\X^3$ and $\gamma_3$ can be illustrated with two simple examples. Consider the double $\dbl\W$ of the ``wedge'' $\W = \{(x,y) \in \R^2~|~ x \geq \abs{y}\}$. $\dbl\W$ is a ``cone'' (see \cref{fig:Wedges}(a)). Suppose that for some small $h > 0$, a geodesic $\gamma : I \to \dbl\W$ starts at the point $(\frac12,h)$ with velocity $(-1,0)$ in one of the copies of $\W$ in $\dbl\W$. Then $\gamma$ will ``wind around'' the vertex of $\dbl\W$ near time $t = \frac12$, and then end at $\gamma(1) = (\frac12, -h)$ (see \cref{fig:Wedges}(a)). During its travel, $\gamma$ crosses into the other copy of $\W$, and then crosses back. If we identify the tangent spaces of $\dbl\W$ at the endpoints of $\gamma$ with $\R^2$, then the parallel transport map of $\gamma$ can be expressed as $-I_2$ where $I_n$ is the $n \times n$ identity matrix (see \cref{fig:Wedges}(b)). In a sense, winding around the vertex of the wedge twists the parallel transport map. Higher-dimensional analogues of this twisting will help us build twisted tubes.

Building upon the previous example, consider the ``3D wedge'' $\R \times \W$ and its double $\dbl(\R \times \W)$ (see \cref{fig:Wedges}(c)). We could consider an analogous geodesic $\tilde\gamma : I \to \dbl(\R \times \W)$ that starts at the point $(0,\frac12,h)$ with velocity $(0,-1,0)$ in one of the copies of $\R \times \W$. The singularities of $\dbl(\R \times \W)$ come from the edge $\R \times \{0\}$ of $\R \times \W$, and $\tilde\gamma$ ``winds around'' this edge near time $t = \frac12$ (see \cref{fig:Wedges}(d)). (When visualizing a geodesic $\tilde\gamma$ in the double of a 3-polytope $\X$, we will adopt the convention of drawing $\X$ and drawing solid line segments for segments of $\tilde\gamma$ in one copy of $\X$ in $\dbl\X$. Dashed segments will represent segments of $\tilde\gamma$ in the other copy of $\X$.) In the coordinates of each copy of $\W$ we may write $\tilde\gamma(t) = (0,\gamma(t))$, so the parallel transport map of $\tilde\gamma$ is $\diag(1,-1,-1)$. More generally, some linear algebra (see \cref{lem:WedgeCollision}) shows that even if the initial velocity of $\tilde\gamma$ is changed to another vector $(v_x,v_y,0)$ such that $v_y < 0$ (so that it travels toward the edge), if we let $\tilde\gamma$ travel until it returns to the copy of $\R \times \W$ that it started from, then its parallel transport map can be expressed as $A \oplus (-1)$, where $A$ is some $2 \times 2$ matrix and $\oplus$ denotes the matrix direct sum.

\begin{figure}[p]
    \centering
    \includegraphics[width=0.7\linewidth]{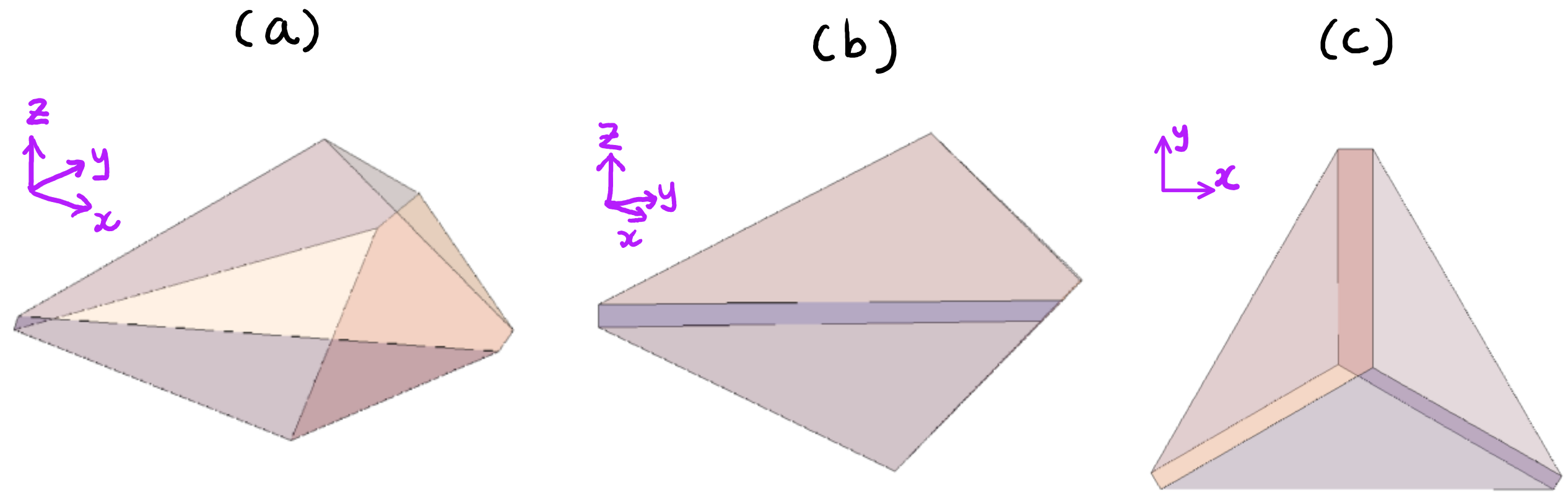}
    \includegraphics[width=0.7\linewidth]{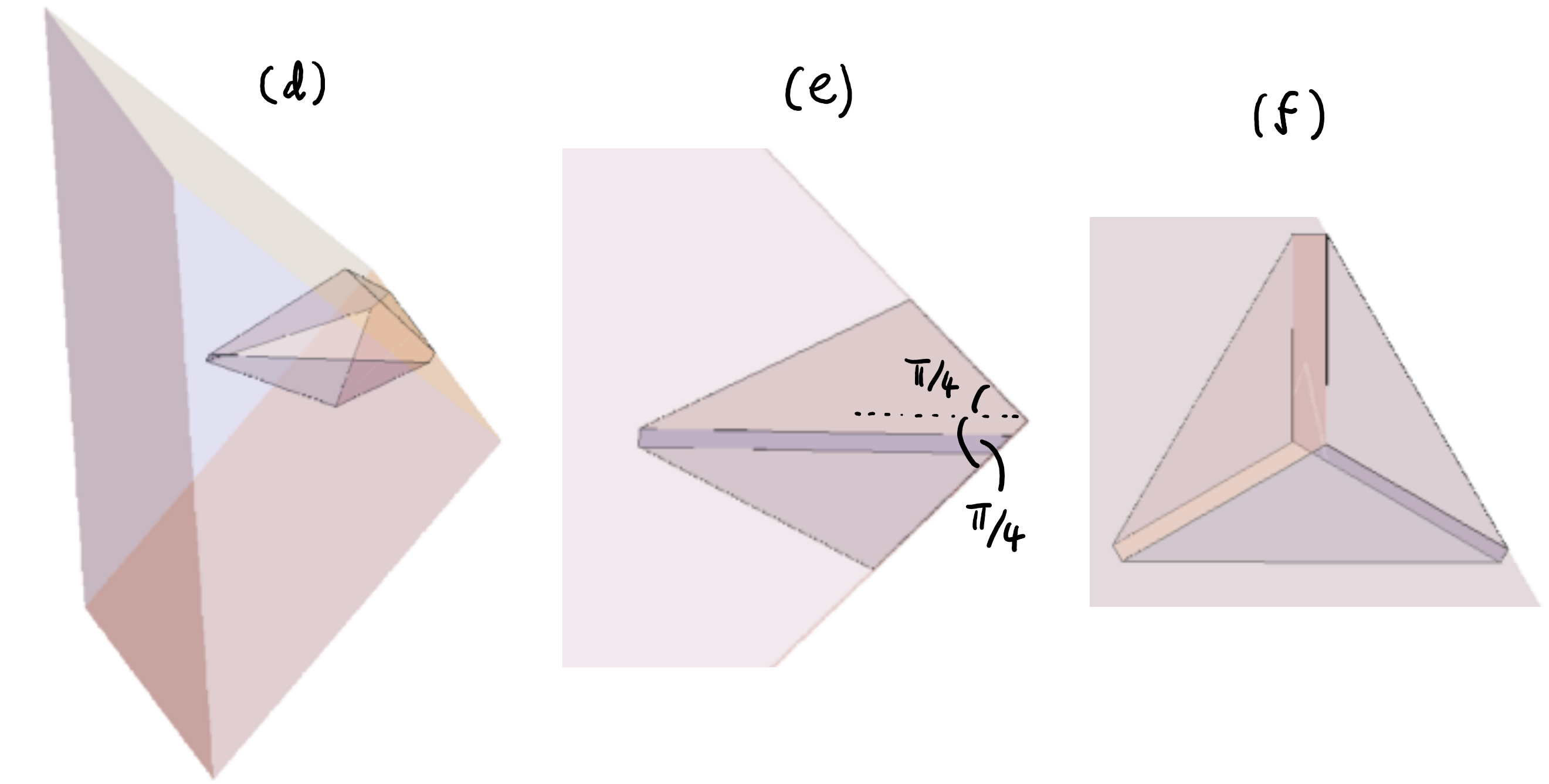}
    \includegraphics[width=0.7\linewidth]{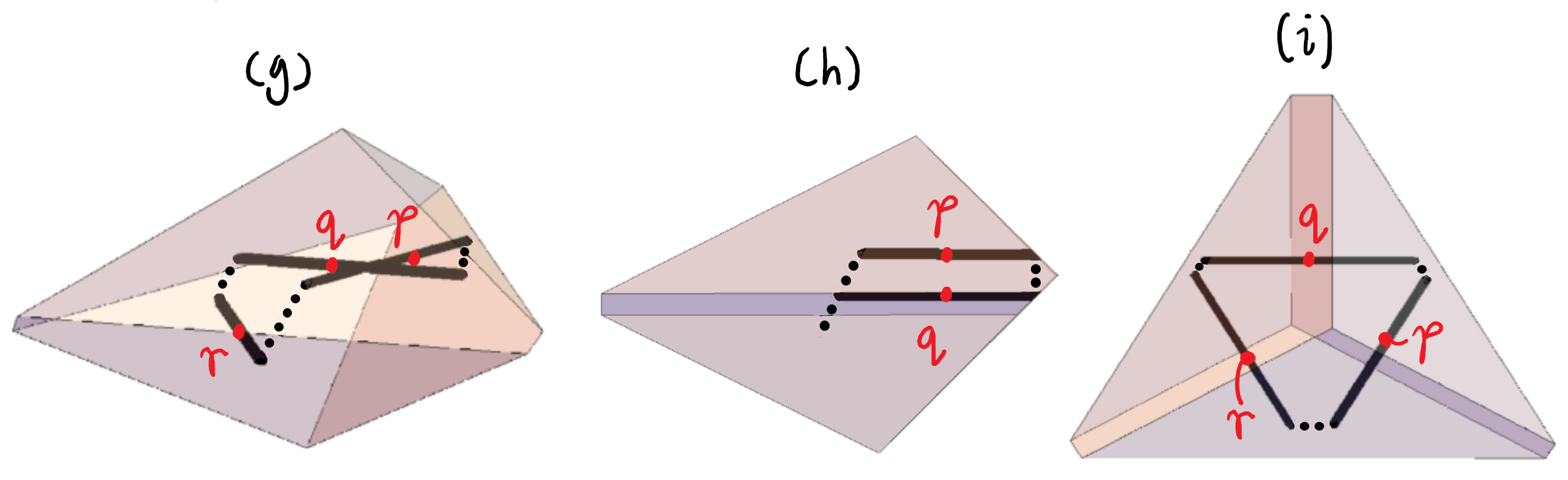}
    \includegraphics[width=0.7\linewidth]{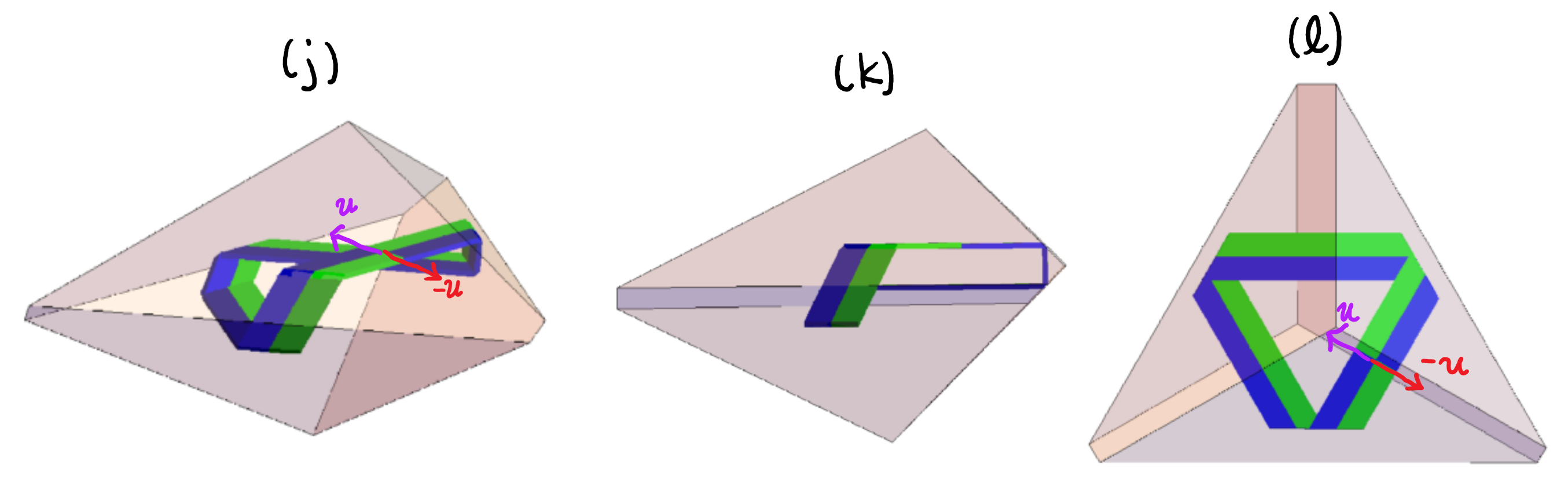}
    \caption{(a)--(c) Three views of $\X^3$. (d)--(f) $\X^3$ is the intersection of three ``wedges'' isometric to $\R \times \W$, and one of the wedges is pictured here. (g)--(i) A closed geodesic $\gamma_3 : I \to \dbl\X^3$. The solid segments lie in one copy of $\X^3$ in $\dbl\X^3$ while the dashed segments lie in the other copy. The solid segments lie parallel to the $xy$-plane. (j)--(l) Parallel transport moves a vector $u$ orthogonal to the geodesic along the green strip, ending at $-u$. The green and blue strips form a M\"obius strip.}
    \label{fig:X3}
\end{figure}

Now we can hint at a proof of the $n = 3$ case of \cref{thm:StableClosedGeodesicPosCurv}: $\X^3$ is the intersection of three isometric copies of $\R \times \W$ (see \cref{fig:X3}(a)--(f); one of the wedges is depicted in (d)--(f)). $\X^3$ is a 3-polytope with six faces, two coming from each wedge. $\gamma_3$ is a closed geodesic in $\dbl\X^3$ (see \cref{fig:X3}(g)--(i)). It can be broken into three piecewise geodesic arcs as shown in \cref{fig:X3}(i): from $p$ to $q$, from $q$ to $r$, and from $r$ to $p$. Each arc collides with the two faces of a single wedge; for instance, the arc from $p$ to $q$ (see \cref{fig:X3}(h)) collides with the faces of the wedge shown in \cref{fig:X3}(e). The previous arguments imply that the arcs have parallel transport maps of the form $A_i \oplus (-1)$ for $2 \times 2$ matrices $A_1$, $A_2$, and $A_3$. Thus $\gamma_3$ itself has parallel transport map $P = A_3A_2A_1 \oplus (-1)^3$. The orientability of $\dbl\X^3$ implies that $\det P = 1$. This, together with the fact that $P$ fixes $\gamma_3'(0)$, guarantees that $P = 1 \oplus (-I_2)$ up to some change of basis. In other words, a tubular neighbourhood of $\gamma_3$ is a twisted tube.

(A more visual way to compute $P$ is to parallel transport a vector $u$ once around $\gamma_3$, where $u$ is orthogonal to $\gamma_3$. \Cref{fig:X3}(j)--(l) shows $u$ being transported along the green strip to $-u$. The green and blue strips form a M\"obius strip, which implies that $P = 1 \oplus (-I_2)$.)

We will prove \cref{thm:StableClosedGeodesicPosCurv}  for $n \geq 4$ inductively, by constructing $\X^n$ as the intersection of $\X^{n-1} \times \R$ with three isometric copies of the ``wedge'' $\R^{n-2} \times \W$. $\gamma_n$ will be carefully chosen so that it ``winds around'' the ``edges'' of those wedges in such a way to yield a parallel transport map that can be expressed in some basis as $1 \oplus (-I_{n-1})$.

\subsection{Stable figure eights from ``incompatibly-twisted'' parallel transport}
\label{sec:Intro_FigureEights}

\Cref{thm:Stable2LoopPosCurv} will be proven by combining the new tool of wedges with earlier techniques from \cite{Cheng_StableGeodesicNets}. Our earlier result in \cite{Cheng_StableGeodesicNets} was proven by controlling the parallel transport maps along the geodesics in a stable geodesic bouquet. Consider a stationary geodesic bouquet $G$ in a flat Riemannian manifold $M$ composed of loops $\gamma_1, \dotsc, \gamma_k : I \to M$ based at a point $p$; since $M$ is flat, a study of the index forms of the $\gamma_i$'s reveals that any variation of $G$ by a vector field will have a non-negative second variation in length. Thus to ensure that $G$ is stable, it suffices to have to prevent the second variation in length from vanishing.

A variation $V$ of $\gamma_i$ (with $V(0) = V(1)$) has vanishing second variation in length only when its component orthogonal to $\gamma$, $V^\perp$, is parallel, because the curvature term vanishes in the formula for the second variation in length. Thus we may associate $\gamma_i$ with a vector space $K_i \subset T_pM$ called its \emph{parallel defect kernel}, which consists of vectors $v$ that can be extended to a variation $V$ of $\gamma_i$ (so $V(0) = V(1) = v$) such that $V^\perp$ is parallel. Then $G$ is stable if $\bigcap_{i = 1}^k K_i = \{0\}$. (This may be viewed as a generalization of twisted tubes to the context of geodesic bouquets.)

In our construction of stable geodesic bouquets in \cite{Cheng_StableGeodesicNets}, the parallel defect kernels turned out to be hyperplanes of dimension $\dim(M) - 1$. Consequently, we needed $\dim(M)$ loops in our geodesic bouquet to guarantee that the intersection of the associated hyperplanes would be a point.

\Cref{thm:Stable2LoopPosCurv} will be proven by constructing $n$-polytopes $\X$ as the intersection of some isometric copies of $\R^{n-2} \times \W$ and some ``prisms'' of the form $\Y \times \R$, where $\Y$ is an $(n-1)$-polytope. These wedges and prisms will be chosen so that $\dbl\X$ will contain geodesic loops whose parallel transport maps give rise to parallel defect kernels have dimension at most $n/2$. The wedges are used to control the parallel transport map more precisely to lower the dimension of the parallel defect kernel. This will eventually allow us to construct stable geodesic bouquets using only two loops by ``twisting'' the parallel transport maps of the two loops in incompatible ways, resulting in parallel defect kernels that intersect in a point. We will carry out the above process for $n = 3,4,5$. Finally, we will combine these low-dimensional constructions into ones of higher dimension using the following technique.

\subsection{Combining low-dimensional constructions into high-dimensional ones}
\label{sec:Intro_Combining}

\begin{figure}[h]
    \centering
    \includegraphics[width=0.9\textwidth]{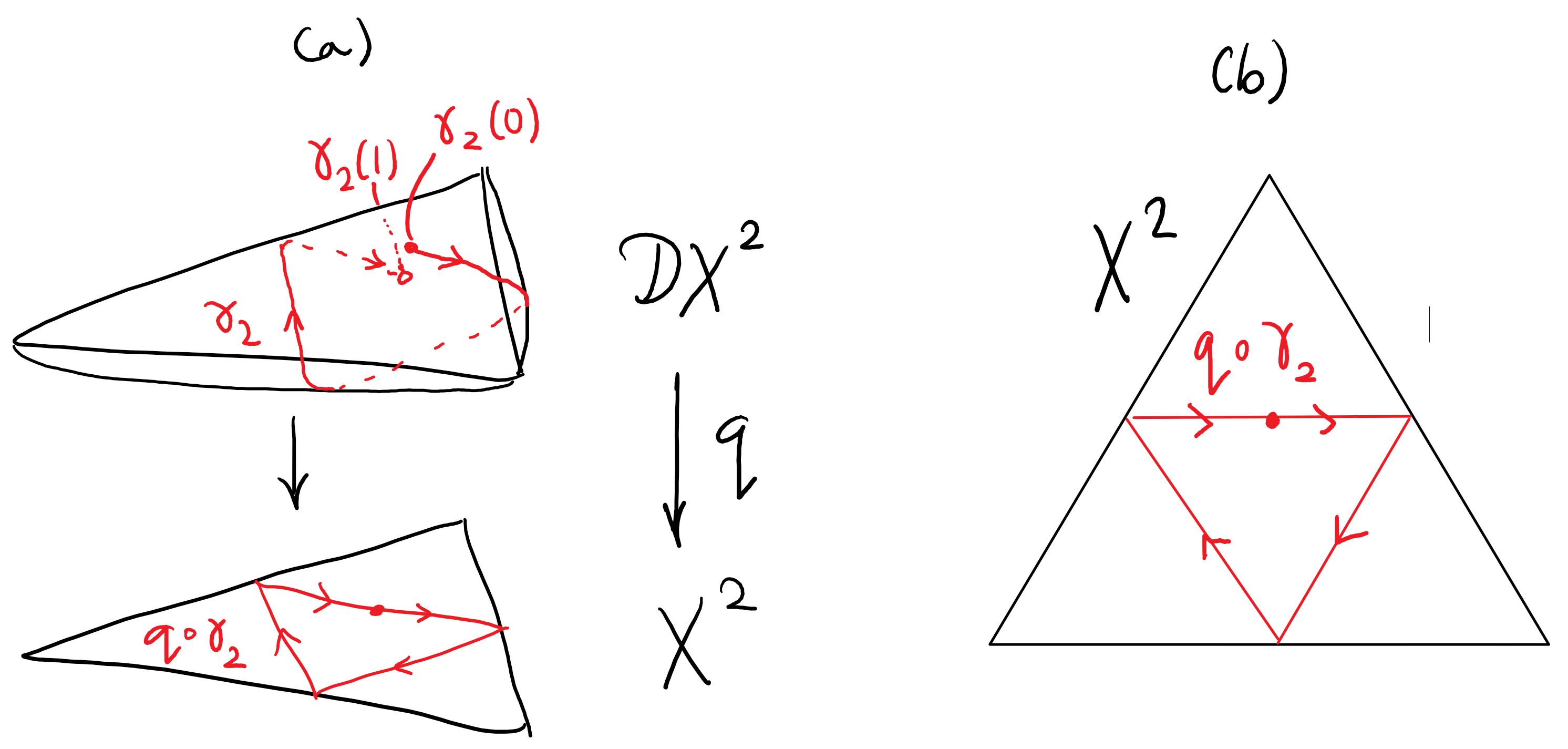}
    \caption{(a) The equilateral triangle $\X^2$ and the geodesic $\gamma_2 : I \to \dbl\X^2$. $q : \dbl\X^2 \to \X^2$ is the canonical quotient map. (b) An illustration of $q \circ \gamma_2$.}
    \label{fig:AntiPeriodicGeodesic}
\end{figure}

If $F$ and $G$ are stable geodesic bouquets with $k$ loops in the doubles of an $m$-polytope $\X_1$ and an $n$-polytope $\X_2$ respectively, then given a few additional minor assumptions we can combine them into a stable geodesic bouquet in the double of an $(m+n)$-polytope. The idea behind this is most naturally presented through the lens of a certain correspondence between billiard trajectories $\beta$ in $\X$ and geodesics $\gamma$ in $\dbl\X$. For instance, it is apparent from \cref{fig:AntiPeriodicGeodesic} that projecting the geodesic $\gamma_2$ via the canonical quotient map $q : \dbl\X^2 \to \X^2$ yields a billiard trajectory in $\X^2$. In general the correspondence is such that $\gamma$ projects to $\beta$ in the same way, and whenever $\beta$ makes a collision, $\gamma$ passes from one copy of $\X$ in $\dbl\X$ to the other.\footnote{The idea behind this correspondence is commonly used in, for example, the study of billiards in rational polygons (i.e. $n = 2$) via their corresponding geodesics in translation surfaces.} (Note that each billiard trajectory could correspond to two possible geodesics, depending on which copy of $\X$ in $\dbl\X$ the geodesic begins from.)

If $\gamma_i : [a,b] \to \dbl\X_i$ are geodesics in the doubles of convex polytopes $\X_i$ for $i \in \{1,2\}$, then they correspond to billiard trajectories $\beta_i : [a,b] \to \X_i$. One can verify that the path $\beta(t) = (\beta_1(t),\beta_2(t))$ is a billiard trajectory in $\X_1 \times \X_2$, \emph{as long as $\beta_1$ and $\beta_2$ never collide simultaneously.} We call $\beta$ a \emph{product billiard trajectory}. Under that assumption, $\beta$ corresponds to some geodesic in $\dbl(\X_1 \times \X_2)$. In this manner we can combine the geodesic loops of $F$ and $G$ into geodesic loops in $\dbl(\X_1 \times \X_2)$, which will form a stationary geodesic bouquet. We will prove that this geodesic bouquet is stable using methods to be introduced later.

\subsection{Organization of content}

\Cref{sec:Definitions} will present definitions important for stationary and stable geodesic bouquets, and establish some notation.

In \cref{sec:BevelingOrigamiModel} we will study how wedges help to twist the parallel transport map of a geodesic. We will use wedges to construct polytopes whose doubles will contain twisted tubes. This will lead to a proof of \cref{thm:StableClosedGeodesicPosCurv}. This section can be read independently from the subsequent sections of the paper.

In \cref{sec:StabilityOfGeodesicBouquets}, we will define parallel defect kernels, relate them to the stability of geodesic bouquets, and compute the parallel defect kernels of the kind of geodesic loops that we will construct to prove \cref{thm:Stable2LoopPosCurv}. In \cref{sec:ProductBilliardTraj} we will analyze the stability of product billiard trajectories using parallel defect kernels.

In \cref{sec:Stable2Loop}, we will construct stable figure eights in the doubles of polytopes with dimensions 3, 4 and 5 using wedges and product billiard trajectories. Finally, we will combine these constructions into stable figure eights in polytopes of all higher dimensions using product billiard trajectories, and complete our proof of \cref{thm:Stable2LoopPosCurv}.

\section{Definitions and Notation}
\label{sec:Definitions}

The definitions in this section originate from \cite[Section~2]{Cheng_StableGeodesicNets}.

\subsection{Definitions for stationary and stable geodesic bouquets}

Let $\bq_k$ denote the ``flower-shaped'' graph with $k$ loops around a single vertex denoted by $*$. A \emph{stationary geodesic bouquet} in a Riemannian manifold $M$ is an immersion $G : \bq_k \to M$ such that each edge is immersed as a non-constant geodesic, and such that the outgoing unit tangent vectors at the \emph{basepoint} $G(*)$ sum to zero.

A \emph{variation} of $G$ is a homotopy $H : (-\varepsilon, \varepsilon) \times \bq_k \to M$ that restricts to a variation of a geodesic on each edge of $\bq_k$. Similarly, a \emph{vector field} along $G$ is a continuous map $V : \bq_k \to TM$ that restricts to a $C^\infty$ vector field along each edge. A vector field $V$ is \emph{tangent to} $G$ if it is tangent to the image of each geodesic loop, and when $k \geq 2$, it also vanishes at the basepoint $G(*)$.\footnote{The reason for this dependence on $k$ is that when $k = 1$, $G$ is a closed geodesic, and stable closed geodesics obey slightly different conditions at the basepoint than for stable geodesic bouquets with $k \geq 2$.}

Adapting a standard definition, for each variation $H$ of $G$, let $\frac{dH}{dt}|_{t = 0}$ denote the vector field along $G$ whose value at $s \in \bq_k$ is the velocity of the curve $t \mapsto H(t,s)$ at time $t = 0$.

A stationary geodesic bouquet $G$ is \emph{stable} when for all variations $H$ of $G$, either $\length_H''(t) > 0$ (where $\length_H(t)$ is the length of the image of $H(t,-) : \bq_k \to M$) or $\frac{dH}{dt}|_{t = 0}$ is tangent to $G$. In this case, $G$ is called a \emph{stable geodesic net}.

\subsection{Notation and other terminology}

We will require that geodesics and billiard trajectories be parametrized at constant speed. We will say that paths $\alpha : [a,b] \to X$ like geodesics or billiard trajectories are \emph{simple} if they are injective except possibly at their endpoints. We may sometimes write $\alpha$ to denote its image. A stationary geodesic bouquet $G : \bq_k \to M$ is called \emph{simple} if $G$ is injective. We will always assume that a billiard trajectory $\beta : [a,b] \to \X$ begins and ends in the interior of $\X$. $\beta$ is called a \emph{billiard loop} if $\beta(a) = \beta(b)$. If $\beta$ also satisfies $\beta'(a) = \beta'(b)$, then it is called a \emph{periodic billiard trajectory}.

When $V : I \to TM$ is a vector field along a geodesic $\gamma : I \to M$, $V(t)^\perp$ will refer to the component of $V(t)$ that is orthogonal to $\gamma'(t)$. $V^\perp$ will refer to the vector field $V^\perp(t) = V(t)^\perp$.

If $f : X \to Y$ and $g : X \to Z$ are maps, let $(f,g)$ denote the map $X \to Y \times Z$ that sends $x$ to $(f(x),g(x))$.

If $F$ is a face of an $n$-polytope, then let $R_F \in O(n)$ denote the reflection about that face. Let $dR_F$ denote the linear part of that reflection.
\section{Stable Closed Geodesics}
\label{sec:BevelingOrigamiModel}

As explained in \cref{sec:Intro_CoreCurves}, to prove \cref{thm:StableClosedGeodesicPosCurv}, for each odd $n \geq 3$, we will find an $n$-polytope $\X^n$ such that $\dbl[sm]\X^n$ contains an isometrically embedded $n$-dimensional twisted tube. Let us prove that its core curve is a stable closed geodesic.

\begin{lemma}
    \label{Lem:MaxTwistStableGeodesic}
    Core curves of twisted tubes are stable closed geodesics.
\end{lemma}
\begin{proof}
    Let $\gamma : [0,\ell] \to T$ be the core curve, parametrized by arc-length, of an $n$-dimensional twisted tube $T$. Consider the formula for the second variation of length \cite[p.~253]{Jost_RiemannianGeometry} along a variational vector field $V(t)$ (with $V(0) = V(\ell)$), noting that it only involves the component of $V$ orthogonal to $\gamma$, and that the curvature term vanishes because $T$ has a flat metric. Thus the second variation is always non-negative, and all we have to prove is that the second variation vanishes only when $V(t)$ is parallel to $\gamma'(t)$ for all $t$.
    
    Let $V(t) = a_1(t)E_1(t) + \dotsb + a_{n-1}(t)E_{n-1}(t) + a_n(t)\gamma'(t)$ for some parallel orthonormal frame $E_1,\dotsc, E_{n-1}, \gamma'$ and smooth functions $a_i$. The definition of a twisted tube $T = \D_r^{n-1} \times [0,\ell]/{\sim}$ implies that we can choose each $E_i(s)$ to be tangent to the image of $\D_r^{n-1} \times \{s\}$ in $T$ and such that $E_i(s)$ corresponds to the $i^\text{th}$ standard basis vector of $\R^{n-1}$ in each copy of $\D_r^{n-1}$. Thus $E_i(\ell) = -E_i(0)$.
    
    If the second variation of length along $V$ vanishes, then
    \begin{equation*}
        0 = \int_0^\ell\norm{\nabla_{\gamma'} \left(\sum_{i = 1}^{n-1} a_i(t)E_i(t)\right)}^2 \,dt = \int_0^\ell\norm{\sum_{i = 1}^{n-1} a_i'(t)E_i(t)}^2 \,dt.
    \end{equation*}
    However, that would imply that the functions $a_1, \dotsc, a_{n-1}$ are all constant. In fact, they must vanish because $V(\ell) = V(0)$ but $E_i(\ell) = -E_i(0)$.
\end{proof}

Finding a twisted tube in a flat manifold is the same as finding a closed geodesic $\gamma$ whose parallel transport map negates all vectors in $\gamma'(0)^\perp$. To study the parallel transport of $\gamma$, it will be convenient to project it down to a billiard trajectory $\beta$ in $\X$, for the following reasons:
\begin{itemize}
    \item The tangent spaces of $\X$ can be canonically identified with $\R^n$.
    
    \item The parallel transport map can be computed in terms of a product of reflections encounted by $\beta$.
\end{itemize}

To formalize this, define the \emph{parallel transport map} of $\beta$ as the linear transformation $P_\beta : T_{\beta(a)}\R^n \to T_{\beta(b)}\R^n$ that is defined by $P_\beta = dR_{F_k} \dotsm dR_{F_1}$, where $\beta$ collides with the faces $F_1, F_2, \dotsc, F_k$ in sequence. Like parallel transport for geodesics, $P_\beta$ is an orthogonal transformation that sends $\beta'(a)$ to $\beta'(b)$. $P_\beta$ is related with the parallel transport map $P_\gamma$ of $\gamma$ by the following identity:

\begin{lemma}
    \label{lem:CompareParallelTransportReflections}
    Let $\X$ be a convex $n$-polytope and $\gamma : [a,b] \to \dbl\X$ be a geodesic such that $\beta = [a,b] \xrightarrow{\gamma} \dbl\X \xrightarrow{q} \X$ is a billiard trajectory. Then,
    \begin{equation}
        dq \circ P_\gamma = P_\beta \circ dq.
    \end{equation}
\end{lemma}

\begin{proof}
    Let us begin with the $k = 1$ case. We may think of $\dbl\X$ as constructed from the space $\X \cup R_{F_1}(\X)$ by gluing $\partial \X$ to $R_{F_1}(\partial\X)$ via the reflection $R_{F_1}$. That is, $\dbl\X \iso (\X \cup R_{F_1}(\X))/{\sim}$, where $x \sim y$ if and only if $x, y \in \partial \X \cup R_{F_1}(\partial\X)$ and $R_{F_1}(x) = y$. Using this construction, the canonical quotient map $q : (\X \cup R_{F_1}(\X))/{\sim} \to \X$ restricts to the identity on $\X$ and restricts to $R_{F_1}$ on $R_{F_1}(\X)$. In this new construction of the double, the geodesic $\gamma$ is a line segment that starts from $\beta(a)$ and travels with constant velocity $\beta'(a)$. Thus parallel transport simply translates vectors along this line segment, and the lemma follows from the new definition of $q$.
    
    The lemma follows for higher values of $k$ by applying the $k = 1$ case to each collision in turn.
\end{proof}

\subsection{Twisting the parallel transport using wedges}

Suppose that $\HS \subset \R^n$ is a half-space. For each $h \in \R$, there are exactly two halfspaces $\HS'$ and $\HS''$ in $\R^{n+1}$ that contain $\HS \times \{h\}$ and whose boundary hyperplanes form an angle of $\pi/4$ with $\HS \times \{h\}$. Then $\HS' \cap \HS''$ is called the \emph{wedge corresponding to $\HS$ at height $h$,} denoted by $\W_h(\HS)$. (This is an isometric copy of $\R^{n-1} \times \W$, where we recall that $\W = \{(x,y) \in \R^2~|~x \geq \abs{y}\}$.)

For each integer $n \geq 1$, let $\pr_n : \R^{n+1} \to \R^n$ denote the projection onto the first $n$ coordinates. The subscript will be dropped when the dimensions are clear from the context.

\begin{lemma}
    \label{lem:WedgeCollision}
    Let $\HS$ be a half-space of $\R^n$ containing a billiard trajectory $\beta : [a,b] \to \HS$ that collides with $\partial\HS$ at time $\frac{a+b}2$. Suppose that a billiard trajectory $\tilde\beta : [a,b] \to \W_h(\HS)$ begins from a point $(\beta(a), h_0)$ for some $h_0 \neq h$, with velocity $(\beta'(a), 0)$. Then $\tilde\beta$ collides once with each face of $\W_h(\HS)$, at times $t_1$ and $t_2$ such that $\frac{t_1 + t_2}2 = \frac{a + b}2$. Moreover, we have the following:
    \begin{enumerate}
        \item For any $t \leq t_1$, $\tilde\beta(t) = (\beta(t), h_0)$.
        
        \item For any $t \geq t_2$, $\tilde\beta(t) = (\beta(t), 2h - h_0)$.

        \item For any $t_1 \leq t \leq t_2$, $\pr(\tilde\beta(t))$ lies at distance $\abs{h - h_0}$ from $\partial\HS$.

        \item $P_{\tilde\beta} = P_\beta \oplus (-1)$.

        \item\label{enum:WedgeCollision_Length} $\beta$ and $\tilde\beta$ have the same length.
    \end{enumerate}
\end{lemma}

\begin{figure}[h]
    \centering
    \includegraphics[width=0.8\linewidth]{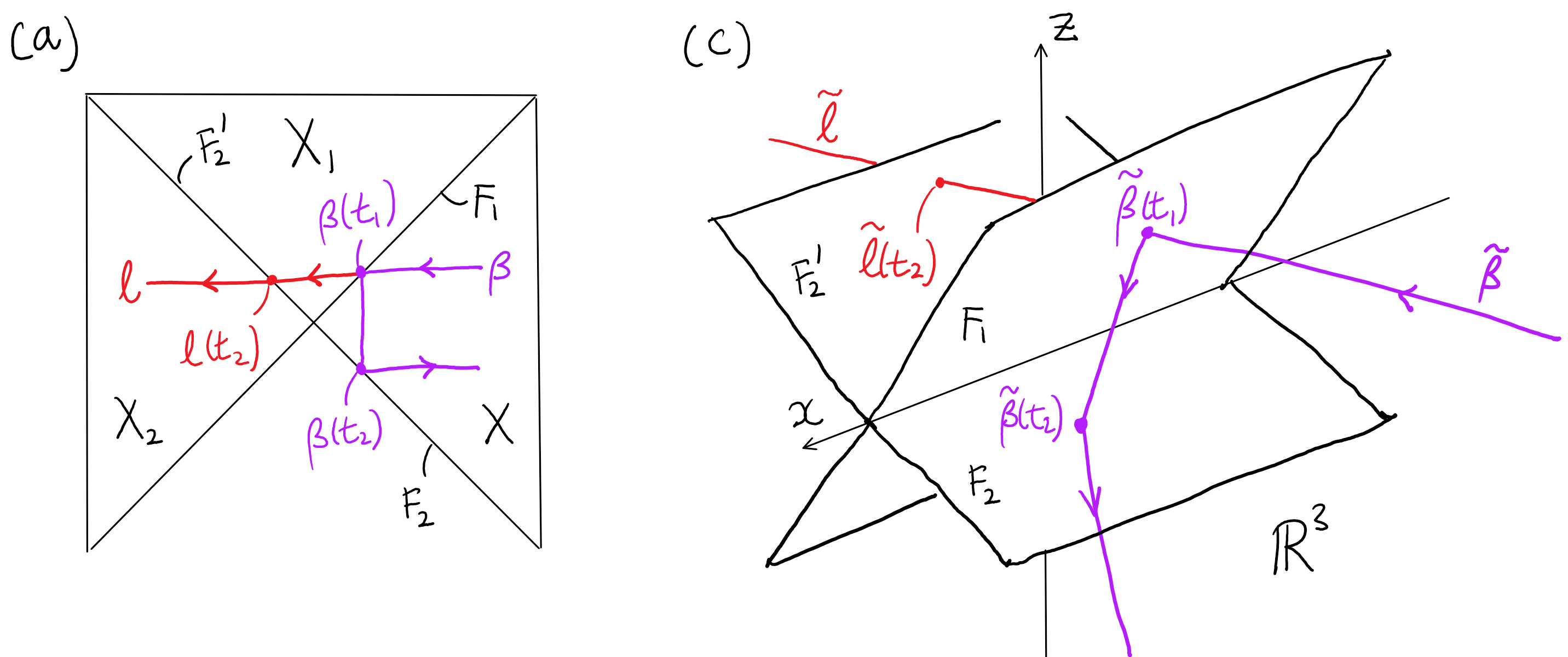}
    \includegraphics[width=0.8\linewidth]{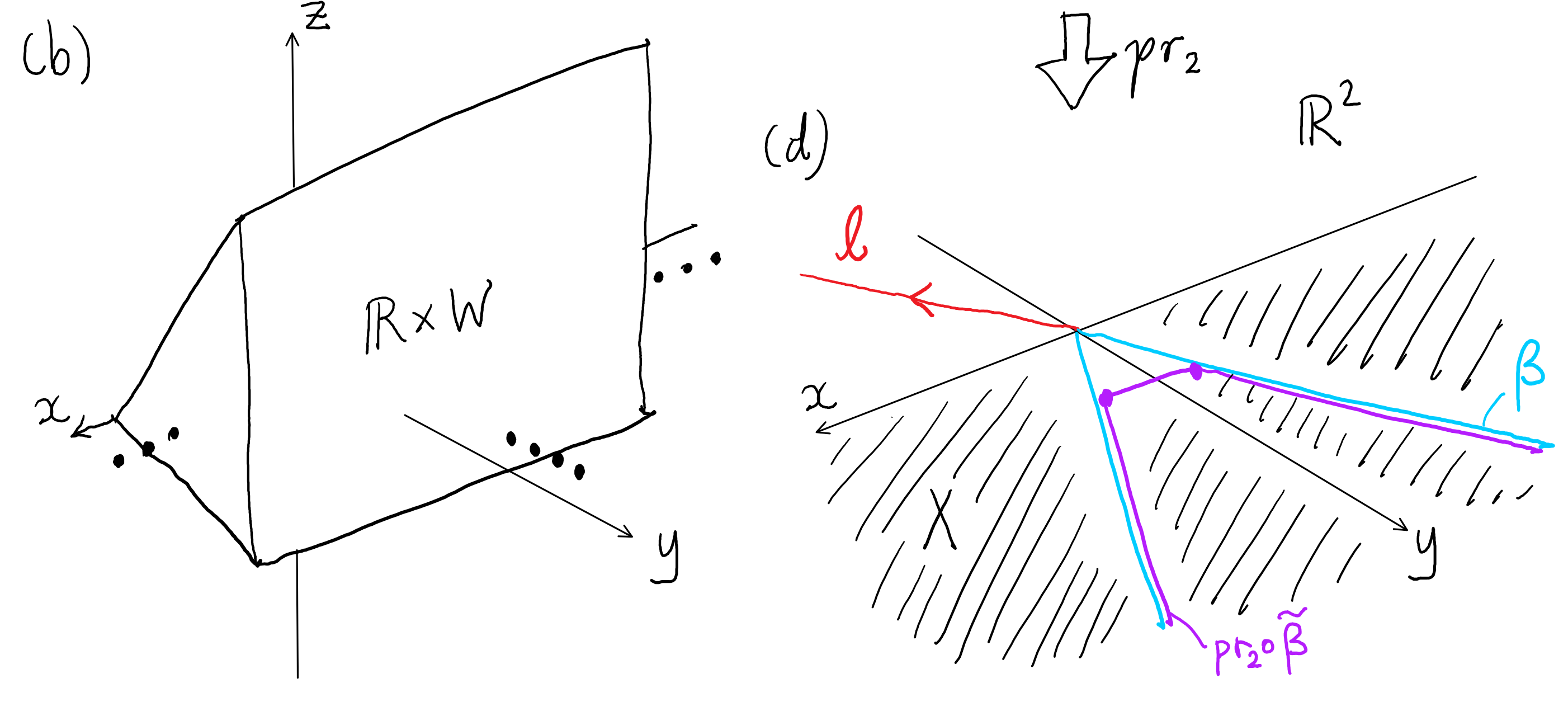}
    \caption{(a) The billiard trajectory $\beta$ in a right-angled isoceles triangle $\X$ unfolds into a line segment $\ell$ (red) that passes through $\X$ and its ``mirrored copies,'' $\X_1$ and $\X_2$. (b) Part of $\R \times \W = \{y \geq \abs{z}\}$, which extends to infinity in three directions as indicated by the ellipses. (c) The billiard trajectory $\tilde\beta$ in $\R \times \W$ unfolds into a line segment $\tilde\ell$ (red). (d) The billiard trajectory $\beta$ (blue) in $\{y \geq 0\}$ unfolds into a line segment $\ell$ (red). The projection $\pr_2 \circ \tilde\beta$ is also shown in purple.}
    \label{fig:Unfolding}
\end{figure}

This lemma can be proven using a tool frequently used to analyze billiard trajectories $\beta : [a,b] \to \X$ in an $n$-polytope $\X$, and sometimes called an \emph{unfolding}. It is inspired by the observation that when bouncing a ball off a mirror, the ball appears to continue into the mirror in a straight line. Suppose that $\beta$ collides at time $t_1 < \dotsb < t_k$ with the faces $F_1, \dotsc, F_k$ of $\X$. We may instead think of $\beta$ as a line segment $\ell : [a,b] \to \R^n$ with the same initial position $\ell(a) = \beta(a)$ and velocity as $\beta$, and which travels through ``mirrored copies'' of $\X$. To be more specific, right after time $t_1$, $\ell$ travels through $\X_1 = R_{F_1}(\X)$ until time $t_2$ when it hits a face $F_2'$ of $\X_1$. Then $\ell$ travels through $\X_2 = R_{F_2'}(\X_1)$ until time $t_3$ when it hits a face $F_3'$ of $\X_2$. This continues until $\ell$ stops at time $b$ in $\X_k$. (See \cref{fig:Unfolding}(a) and (c).)

Some elementary geometry shows that $\X_i = R_{F_1}R_{F_2}\dotsm R_{F_i}(\X)$, which allows us to construct $\beta$ by reflecting portions of $\ell$:
\begin{equation}
    \label{eq:Unfolding}
    \beta(t) = \begin{cases}
        \ell(t) & a \leq t \leq t_1
        \\
        R_{F_1}(\ell(t)) & t_1 \leq t \leq t_2
        \\
        R_{F_2}R_{F_1}(\ell(t)) & t_2 \leq t \leq t_3
        \\
        \qquad\quad \vdots & \qquad\, \vdots
        \\
        R_{F_k}\dotsm R_{F_1}(\ell(t)) & t_k \leq t \leq b.
    \end{cases}
\end{equation}

\begin{proof}[Proof of \cref{lem:WedgeCollision}]
    By isometry, it suffices to prove this lemma in the case that $\HS = \{x_n \geq 0\}$, $h = 0$, and $\W_h(\HS) = \R^{n-2} \times \W$. By a reflection $x_{n+1} \mapsto -x_{n+1}$ in $\R^{n+1}$, we may also assume that $h_0 > 0$. By reparametrizing $\beta$ at a different speed, we may also assume that $[a,b] = I$. (This is the situation pictured in \cref{fig:Unfolding}(b)--(d).)
    
    Consider the unfolding of $\beta$ into a line segment $\ell : I \to \R^n$, which passes through $\X$ and its mirrored copy, $\{x_n \leq 0\}$ (see \cref{fig:Unfolding}(d)). The unfolding of $\tilde\beta$ is a line segment $\tilde\ell : I \to \R^{n+1}$ that passes through $\R^{n-2} \times \W$, then its mirrored copy $\W_1 = \{x_{n+1} \geq \abs{x_n}\}$, and finally yet another copy $\W_2 = \{x_n \leq -\abs{x_{n+1}}\}$ (see \cref{fig:Unfolding}(c)). Indeed, the fact that $\tilde\ell$ starts in the interior of $\R^{n-2} \times \W$ means that its $x_n$-coordinate starts off at a value $r > h_0$. Then its $x_n$-coordinate decreases linearly to $-r$ (this is due to symmetry: the $x_n$-coordinate of $\beta$ reaches 0 when it collides at time $\frac{a+b}2$) while its $x_{n+1}$-coordinate stays constant at $h_0$, which is why $\tilde\ell$ must pass through $\W_1$ and then $\W_2$. This also implies that $\tilde\ell$ intersects the hyperplanes $\Pi_+ = \{x_n = x_{n+1}\}$ once at $t_1$ and $\Pi_- = \{x_n = -x_{n+1}\}$ once at $t_2$, where $\frac{t_1 + t_2}2 = \frac{a + b}2$.
    
    Since $h_0 > 0$, $\tilde\beta$ collides with the faces of $\R^{n-2} \times \W$ exactly once each, first $F_1 \subset \Pi_+$ and then $F_2 \subset \Pi_-$. Therefore we may see that $R_{F_1} = I_{n-2} \oplus \left[\begin{smallmatrix}
        0 & 1\\
        1 & 0
    \end{smallmatrix}\right]$ and $R_{F_2} = I_{n-2} \oplus \left[\begin{smallmatrix}
        0 & -1\\
        -1 & 0
    \end{smallmatrix}\right]$. This allows us to prove (4):
    \begin{equation*}
        P_{\tilde\beta} = dR_{F_2}dR_{F_1} = R_{F_2}R_{F_1} = I_{n-2} \oplus (-I_2) = R_{\partial\HS} \oplus (-1) = dR_{\partial\HS} \oplus (-1) = P_\beta \oplus (-1).
    \end{equation*}
    This calculation also yields (2) via an application of \cref{eq:Unfolding}: when $t \geq t_2$,
    \begin{equation*}
        \tilde\beta(t) = R_{F_2}R_{F_1}(\tilde\ell(t)) = \left(R_{\partial\HS} \oplus (-1)\right)(\ell(t),h_0) = (R_{\partial\HS}(\ell(t)),-h_0) = (\beta(t), -h_0).
    \end{equation*}
    
    \Cref{fig:Unfolding}(c)--(d) show intuitively that (1)--(3) hold. The definition of $\tilde\beta$ implies that $\tilde\ell(t) = (\ell(t),h_0)$ for all $t \in I$, from which (1) follows. The fact that $\tilde\beta$ is derived from $\tilde\ell$ by applying some reflections to certain segments implies that $\tilde\beta$ is as long as $\tilde\ell$, which is as long as $\beta$, giving (\ref{enum:WedgeCollision_Length}).
    
    Finally, (3) follows from the fact that $\tilde\ell$ is at a constant distance of $h_0$ from the hyperplane $\{x_{n+1} = 0\}$. Thus \cref{eq:Unfolding} implies that when $t_1 \leq t \leq t_2$, $\tilde\beta(t) = R_{F_1}(\tilde\ell(t))$ lies at a constant distance of $h_0$ from the hyperplane $\{x_n = 0\} = R_{F_1}(\{x_{n+1} = 0\})$.
\end{proof}

The following proposition will allow us to construct periodic billiard trajectories in convex polytopes of arbitrarily high dimension while controlling their parallel transport maps.

\begin{proposition}
    \label{prop:Beveling}
    Let $\X$ be an $n$-polytope containing a billiard trajectory $\beta : I \to \X$ that collides at most once with each face. Choose any set of $k \geq 2$ faces that $\beta$ collides with, and let the half-spaces containing $\X$ and bounded by those faces be $\HS_1,\dotsc,\HS_k$. Then there exists some $h_1, \dotsc, h_k \in \R$ such that the $(n+1)$-polytope
    \begin{equation}
        \label{eq:Beveling}
        \tilde\X = \X \times \R \cap \W_{h_1}(\HS_1) \cap \dotsb \cap \W_{h_k}(\HS_k)
    \end{equation}
    contains a billiard trajectory $\tilde\beta : I \to \tilde\X$ that satisfies the following properties:
    \begin{enumerate}
        \item $\tilde\beta(t) = (\beta(t),0)$ and $\tilde\beta'(t) = (\beta'(t),0)$ for $t = 0, 1$.

        \item $\tilde\beta$ makes $k$ more collisions than $\beta$.

        \item $\tilde\beta$ collides with each face of $\tilde\X$ at most once.

        \item $P_{\tilde\beta} = P_\beta \oplus (-1)^k$.

        \item If $\beta$ is simple then $\tilde\beta$ is also simple.

        \item Every collision time of $\tilde\beta$ is within distance $\tau > 0$ from a collision time of $\beta$, where $\tau$ can be made arbitrarily small. Moreover, the collision times of $\tilde\beta$ with the boundaries of wedges $\W_{h_i}(\HS_i)$ can be arranged to avoid any finite set of real numbers.

        \item\label{enum:Beveling_Length} $\beta$ and $\tilde\beta$ have the same length.
    \end{enumerate}
\end{proposition}
\begin{proof}
    Without loss of generality, let us order the $\HS_k$'s in the order that their boundaries are encountered by $\beta$. Pick some real numbers $h_1',\dotsc, h_k'$, and let $\varepsilon > 0$ be a small parameter whose magnitude will be chosen later. Define $h_i = \varepsilon h_i'$, $\eta_i = 2(h_i - h_{i-1} + \dotsb + (-1)^{i-1}h_1)$, and $\eta_0 = 0$. Choose the values of $h_i'$ so that $\eta_k = 0$; this is possible because $k \geq 2$. We can now define $\tilde\X$ by \cref{eq:Beveling}. The sizes of the $\eta_i$'s are bounded by the parameter $\eta^* = 2\varepsilon(\abs{h_1'} + \dotsb + \abs{h_k'})$.
    
    \Cref{fig:X3}(a)--(c) shows what $\tilde\X$ could look like when $n = k = 3$. (One of the wedges $\W_{h_i}(\HS_i)$ is shown in \cref{fig:X3}(d)--(f).)
    Choose a disjoint set of time intervals $[a_i, b_i] \subset I$ for $i = 1,\dotsc, k$, such that $\beta$ collides with $\partial\HS_i$ at time $\frac{a_i + b_i}2$ and collides with no other face of $\X$ during the time interval $[a_i, b_i]$. In addition, choose $a_i > 0$ and $b_k < 1$.

    Let us begin by defining $\tilde\beta$ over the times $J = I \setminus \bigcup_i(a_i,b_i)$. For each $i = 0, \dotsc, k$ and $b_i \leq t \leq a_{i+1}$, define $\tilde\beta(t) = (\beta(t),\eta_i)$. (We define $b_0 = 0$ and $a_{k+1} = 1$.) Since $\eta_k = 0$, this satisfies the property (1) in the statement of the proposition. We must ensure that $\tilde\beta(J)$ lies in $\tilde\X$; since it clearly lies in $\X \times \R$ by construction, it lies in $\tilde\X$ as a result of the following lemma:
    
    \begin{lemma}
        \label{lem:AvoidingWedges}
        For each $i = 1,\dotsc,k$, let $J_i = I \setminus (a_i, b_i)$ and consider the convex hull $K$ of $\beta(J_i)$. Then for sufficiently small $\varepsilon$, $K \times [-\eta^*, \eta^*]$ lies in the interior of $\W_{h_j}(\HS_j)$ for all $j \neq i$.
    \end{lemma}
    \begin{proof}
        The assumption that $\beta$ only collides with each face of $\X$ at most once, and that $\X$ is convex, implies that $\beta(J_i)$ lies at some positive distance from $\partial\HS_i$. The same is true of $K$, so $K \times \{0\}$ lies in the interior of $\W_0(\HS_j)$. When $\varepsilon$ (and thus $h_j$ and $\eta^*$) is sufficiently small, $K \times [-\eta^*, \eta^*]$ and $\W_{h_j}(\HS_j)$ will be close enough to $K \times \{0\}$ and $\W_0(\HS_j)$ that $K \times [-\eta^*, \eta^*]$ will lie in the interior of $\W_{h_j}(\HS_j)$.
    \end{proof}

    \Cref{lem:AvoidingWedges} implies that for each $i = 1,\dotsc,k$, $\tilde\beta(a_i) = (\beta(a_i), \eta_{i-1})$ lies in the interior of $\W_{h_i}(\HS_i)$, and $\tilde\beta'(a_i) = (\beta'(a_i),0)$. Since $\beta$ collides with $\partial\HS_i$ at time $\frac{a_i+b_i}2$, \cref{lem:WedgeCollision} allows us to extend $\tilde\beta$ continuously over $[a_i,b_i]$ so that $\tilde\beta|_{[a_i,b_i]}$ is a billiard trajectory in $\W_{h_i}(\HS_i)$, and the velocity $\tilde\beta'(t)$ is continuous at times $a_i$ and $b_i$. (Here we use \cref{lem:WedgeCollision}(2) and the fact that $2h_i - \eta_{i-1} = \eta_i$.)

    Thus we have defined $\tilde\beta$ over $I$, and we need to check that $\tilde\beta(I) \subset \tilde\X$. For each $1 \leq i, j \leq k$ such that $i \neq j$, observe that $\tilde\beta([a_i,b_i])$ lies in the convex hull of $\beta(J \cup [a_i,b_i]) \times [-\eta^*, \eta^*]$. This implies that $\tilde\beta([a_i,b_i]) \subset \X \times \R$, but \cref{lem:AvoidingWedges} also guarantees that for sufficiently small $\varepsilon$, $\tilde\beta([a_i,b_i])$ lies in the interior of $\W_{h_j}(\HS_j)$. Therefore $\tilde\beta(I) \subset \tilde\X$.

    Now we verify that $\tilde\beta$ is a billiard trajectory in $\tilde\X$. $\tilde\beta|_J$ touches $\partial(\X \times \R)$ at times matching the collision times of $\beta|_J$. The touching points lie in the interiors of the wedges $\W_{h_i}(\HS_i)$, due to \cref{lem:AvoidingWedges}. \Cref{lem:WedgeCollision}(3) implies that for each $i = 1,\dotsc,k$, $\tilde\beta([a_i,b_i])$ lies in the interior of $\X \times \R$. As explained earlier, $\tilde\beta([a_i,b_i])$ lies in the interiors of the wedges $\W_{h_j}(\HS_j)$ for $j \neq i$. Altogether, this implies that $\tilde\beta$ touches the faces of $\tilde\X$ only in their interiors. The rest is straightforward to verify that $\tilde\beta$ is a billiard trajectory in $\tilde\X$. Moreover, the collisions of $\tilde\beta$ consist of two collisions with each $\partial\W_{h_i}(\HS_i)$, and one for each collision of $\beta|_J$. This yields property (2).
    
    By choosing new values of $a_i$ and $b_i$ to shrink the intervals $[a_i, b_i]$ while fixing their midpoints---which requires shrinking $\varepsilon$ even further---we can ensure that the collisions of $\tilde\beta$ happen at times within distance $\tau > 0$ to the collision times of $\beta$, where $\tau$ can be made arbitrarily small. We know from \cref{lem:WedgeCollision} that the collision times of $\tilde\beta$ avoid the real numbers $\frac{a_i + b_i}2$ for all $i$. Given a finite set $S$ of other real numbers, the collision times of $\tilde\beta$ with the boundaries of wedges can be arranged to avoid $S$ by sufficiently shrinking $\tau$. This proves property (6).

    To prove property (5), note that $\tilde\beta$ consists of ``horizontal arcs'' $\tilde\beta|_{[b_i,a_{i+1}]}$ whose $x_{n+1}$-coordinates remain constant at $\eta_i$, and ``transitional arcs'' $\tilde\beta|_{[a_i,b_i]}$ whose $x_{n+1}$-coordinates change from $\eta_{i-i}$ to $\eta_i$. Each horizontal and transitional arc is simple, and moreover the horizontal arcs cannot intersect each other (except that the first and last one will intersect at $\tilde\beta(0)$) because they project via $\pr$ to arcs of $\beta$, whereas $\beta$ is simple. By shrinking the intervals $[a_i, b_i]$ while fixing their midpoints as before, we can arrange that each transitional arc intersects no other arc except for the horizontal arcs right before and after it.

    Next, we prove property (3). For each time $t \in J$ when $\beta$ collides with a face $F$ of $\X$, $\tilde\beta$ collides with a face of $\tilde\X$ that is contained in $F \times \R$. The segments $\tilde\beta|_{[a_i,b_i]}$ also collide with the faces of $\W_{h_i}(\HS_i)$ once each. Thus $\tilde\beta$ never collides with a given face more than once.

    Property (\ref{enum:Beveling_Length}) follows from observing that $\tilde\beta|_J$ and $\beta|_J$ have the same length, and that each $\tilde\beta|_{[a_i,b_i]}$ has the same length as $\beta|_{[a_i,b_i]}$ by \cref{lem:WedgeCollision}(\ref{enum:WedgeCollision_Length}).

    Finally, let us prove property (4). For any $0 \leq s < t \leq 1$, let $P_{st}$ denote the parallel transport map of $\beta|_{[s,t]}$, if defined. For any $i = 0,\dotsc,k$, $\tilde\beta|_{[b_i, a_{i+1}]}$ has parallel transport map $P_{b_ia_{i+1}} \oplus 1$. For any $i = 1,\dotsc,k$, $\tilde\beta|_{[a_i, b_i]}$ has parallel transport map $P_{a_ib_i} \oplus (-1)$ by \cref{lem:WedgeCollision}(4). Multiplying them together,
    \begin{equation*}
        P_{\tilde\beta} = P_{b_k1}P_{a_kb_k}P_{b_{k-1}a_k}\dotsm P_{a_1b_1}P_{0a_1} \oplus (-1)^k = P_\beta \oplus (-1)^k.
    \end{equation*}
\end{proof}

\subsection{The construction of stable closed geodesics}
\label{sec:StableClosedGeodesics}

We are almost ready to prove \cref{thm:StableClosedGeodesicPosCurv}, except that we will need the following technical results from our prior work; they guarantee that the double of a convex polytope containing an stable closed geodesic can be smoothed into a convex hypersurface while preserving the stability of that closed geodesic.

\begin{proposition}[{\cite[Appendix~C]{Cheng_StableGeodesicNets}}]
    \label{prop:Smoothing}
    Let $\X$ be a convex $n$-polytope such that $\dbl\X$ contains a stable geodesic bouquet $G$. Then for any compact neighbourhood $N$ of $G$ in $\dbl[sm]\X$, there exists a sequence of embeddings $\{\varphi_i : N \to M_i\}_{i = 1}^\infty$ into smooth convex hypersurfaces $M_i$ of $\R^{n+1}$ with strictly positive curvature, such that the pullback metrics $g_i$ along $\varphi_i$ from $M_i$ converge to the flat metric on $N$ in the $C^\infty$ topology.
\end{proposition}

\Cref{prop:Smoothing} was proved using convolutions, which are widely used to smooth convex bodies and convex functions \cite{Schneider_SmoothApproxConvex,Schmuckenschlaeger_ApproxConvex,Ghomi_SmoothingConvex}.

\begin{proposition}[{\cite[Appendix~B]{Cheng_StableGeodesicNets}}]
    \label{prop:TransferGeodesicBouquet}
    Let $(N,g_0)$ be a compact and flat Riemannian manifold, whose interior contains a simple, irreducible, and stable geodesic bouquet $G_0$. Then for any Riemannian metric $g$ on $N$ that is sufficiently close to $g_0$ in the $C^\infty$ topology, $(N,g)$ also contains an irreducible and stable geodesic bouquet with the same number of loops.
\end{proposition}

\Cref{prop:TransferGeodesicBouquet} was proved by adapting well-known methods to approximate path spaces by finite-dimensional manifolds \cite{Milnor_MorseTheory}, as well as results showing that the critical points of certain Morse functions persist when the function is slightly perturbed  \cite{GolubitskyGuillemin_StableMappings,Ruas_OldNewResultStableMappings}. 

Now we are ready to prove \cref{thm:StableClosedGeodesicPosCurv}.

\begin{proof}[Proof of \cref{thm:StableClosedGeodesicPosCurv}]
    We begin by inductively constructing, for every integer $n \geq 2$, an $n$-polytope $X^n$ that has a simple and periodic billiard trajectory $\beta_n : I \to \X^n$ whose parallel transport map is $1 \oplus (-I_{n-1})$ in some basis. Moreover, $\beta_n$ collides $3(n-1)$ times, and it collides with each face of $\X^n$ at most once.

    Let $\X^2$ be an equilateral triangle, and let $\beta_2$ be a periodic billiard trajectory that collides once at each midpoint of an edge of $\X^2$ (see \cref{fig:AntiPeriodicGeodesic}). Multiplying the reflections about each edge shows that $\beta_2$ has the required parallel transport map.

    Once we have constructed $\X^n$ and $\beta_n$, apply \cref{prop:Beveling} to any choice of 3 faces of $\X^n$ that $\beta_n$ collides with to construct $\X^{n+1}$ and a billiard trajectory $\beta_{n+1} : I \to \X^{n+1}$ with the required properties.

    For each odd integer $n \geq 3$, $\beta_n$ has an even number of collisions, and is a simple closed curve. Hence it corresponds to a simple closed geodesic $\gamma$ in $\dbl\X^n$. The parallel transport map is $1 \oplus (-I_{n - 1})$ with respect to the decomposition of $\R^n$ given by a tangent vector of $\beta_n$ and its orthogonal complement, which implies that a tubular neighbourhood of $\gamma$ is in fact an $n$-dimensional twisted tube $T$. Hence, \cref{Lem:MaxTwistStableGeodesic} implies that $\gamma$ is stable.
    
    Since $\gamma$ is simple and it is also a stable geodesic bouquet with one loop, we may apply \cref{prop:Smoothing,prop:TransferGeodesicBouquet} to prove the desired result.
\end{proof}
\section{Stable figure eights}

\subsection{Sufficient conditions for stability}
\label{sec:StabilityOfGeodesicBouquets}

The stability of closed geodesics can be deduced by studying their bilinear index forms $I_\gamma(-,-)$ \cite[Definition~6.1.1]{Jost_RiemannianGeometry}. We will generalize this to stationary geodesic bouquets $G$ by considering a quadratic form $Q_G(-)$, also called its \emph{index form}, which we define as follows.\footnote{This definition is equivalent to the one made previously in \cite[Equation~(3.2)]{Cheng_StableGeodesicNets}.} if $V$ is a vector field along $G$ that restricts to $V_i$ along $\gamma_i$, then
\begin{equation}
    \label{eq:IndexFormDef_GeodesicBouquet}
    Q_G(V) = \sum_{i = 1}^k I_{\gamma_i}(V_i^\perp,V_i^\perp) = \sum_{i=1}^k \int_0^1 \norm{\nabla_{\gamma_i'}V_i^\perp}^2 - \ip{R(\gamma_i',V_i^\perp)V_i^\perp,\gamma_i'} \,dt,
\end{equation}
where $R$ is the Riemann curvature tensor. We will abbreviate $I_\gamma(V,V)$ as $Q_\gamma(V)$.

Index forms are related to stability by the following lemma that generalizes the same result for closed geodesics. It is a direct consequence of the definitions.
\begin{lemma}[{\cite[p.~9]{Cheng_StableGeodesicNets}}]
    \label{lem:IndexFormSecondVariation}
    Let $G$ be a stationary geodesic bouquet in a flat Riemannian manifold $M$. Let $V$ be a vector field along $G$. Then any variation $H$ of $G$ in the direction of $V$ satisfies $\length_H''(0) \geq 0$, where
    \begin{equation}
        \label{eq:2ndVariationGeodesicBouquet}
        \length_H''(0) = 0 \iff Q_G(V) = 0.
    \end{equation}
    Moreover, $G$ is stable if and only if the nullspace of $Q_G$ contains only vector fields tangent to $G$.
\end{lemma}

To prove that a stationary geodesic bouquet is stable, we cannot apply our knowledge of maximally-twisted tubes; we will generalize this theory to find conditions on parallel transport along the geodesics of a stationary figure eight that would guarantee stability. Such conditions were identified in \cite[Section~3]{Cheng_StableGeodesicNets}. In this section, we will briefly recapitulate the portion of this theory that is essential for us. We will include their proofs to keep the exposition self-contained. We will then apply this theory to describe the parallel defect kernels of geodesic loops of the kind that we will construct in that proof.

As explained in \cref{sec:Intro_FigureEights}, each geodesic loop $\gamma$ in a stationary geodesic bouquet can be associated with a subspace of the tangent space at the basepoint called the \emph{parallel defect kernel}. It is the kernel of a \emph{parallel defect operator}, defined as follows.

\begin{definition}[Parallel defect operator {\cite{Cheng_StableGeodesicNets}}]
    \label{def:ParallelDefectOperator}
    Let $\gamma : [0,1] \to M$ be a geodesic loop. For $t = 0,1$, let $\pi_t : T_{\gamma(t)}M \to T_{\gamma(t)}M$ be the projection onto $\gamma'(t)^\perp$. Then the \emph{parallel defect operator of $\gamma$} is $\opd\gamma = P_\gamma\pi_0 - \pi_1$.
\end{definition}

In some sense, the quantity $\norm{\opd\gamma v}$ measures how difficult it is to extend a single vector $v \in T_pM$ to a vector field $V$ along $\gamma$ so that $V^\perp$ is parallel along $\gamma$---that is, so that $Q_\gamma(V^\perp) = 0$. This heuristic is formalized in the following inequality.

\begin{proposition}[{\cite[p.~10]{Cheng_StableGeodesicNets}}]
    \label{prop:IndexFormLowerBound_ParallelDefect_Flat}
    Let $M$ be a flat Riemannian manifold. Suppose that $M$ contains a geodesic loop $\gamma : [0,1] \to M$. Then for any vector field $V$ along $\gamma$,
    \begin{equation}
        \label{eq:IndexFormLowerBound_GeodesicLoop_Flat}
        Q_\gamma(V^\perp) \geq \norm{\opd\gamma V(0)}^2.
    \end{equation}
\end{proposition}
\begin{proof}
    Choose a parallel orthonormal frame $F_1, \dotsc, F_{n-1},\gamma'$ along $\gamma$, where $n = \dim M$. Let $V^\perp = \sum_{i = 1}^{n-1}c_iF_i$ for piecewise smooth functions $c_i : [0,1] \to \R$. Observe that $Q_\gamma(V^\perp) = \int_0^1 \sum_{i = 1}^{n-1}c_i'(t)^2\,dt$ is the twice the energy of the path $t \mapsto (c_1(t),\dotsc, c_{n-1}(t))$ in $\R^{n-1}$ \cite[eq.~(1.4.8)]{Jost_RiemannianGeometry}. If we fix the endpoints of this path, then its energy is minimized by the constant-speed straight path between those endpoints. This minimum energy is $\frac12\sum_{i = 1}^{n-1} (c_i(1) - c_i(0))^2$, which is equal to $\frac12\norm{\opd\gamma V(0)}^2$; this follows from unravelling the definition of $\opd\gamma V(0)$ and noting that $V(0) = V(1)$. Thus we obtain \cref{eq:IndexFormLowerBound_GeodesicLoop_Flat}.
\end{proof}

This gives the following sufficient condition for stability.\footnote{In a flat manifold, it is also a necessary condition, but we will not prove that as we will not need it.}

\begin{corollary}[{\cite[p.~11]{Cheng_StableGeodesicNets}}]
    \label{cor:NullvarsToIntersectKernels_Flat}
    Let $G$ be a stationary geodesic bouquet in a flat manifold, based at $p$, with loops $\gamma_1, \dotsc, \gamma_k$. Then $G$ is stable if $\bigcap_{i = 1}^k \ker\opd{\gamma_i} = \{0\}$.
\end{corollary}
\begin{proof}
    To show that $G$ is stable it suffices to check that the nullspace of $Q_G$ contains only vector fields tangent to $G$, by virtue of \cref{lem:IndexFormSecondVariation}. Suppose now that $\bigcap_i\opv{\gamma_i} = \{0\}$, and that $Q_G(V) = 0$. Since the ambient manifold is flat, each $Q_{\gamma_i}$ is non-negative definite, therefore we get $Q_{\gamma_i}(V|_{\gamma_i}^\perp) = 0$. Let $V$ take the value of $v$ at the basepoint of $G$; then \cref{prop:IndexFormLowerBound_ParallelDefect_Flat} implies that $\opd{\gamma_i}v = 0$, and this holds for every $i$. By hypothesis, we have $v = 0$. This means that $V|_{\gamma_i}^\perp$ vanishes at the basepoint, but it is parallel (because $Q_{\gamma_i}(V|_{\gamma_i}^\perp) = 0$) so it must vanish along the entire $\gamma_i$. Since the above holds for all $i$, $V$ must be tangent to $G$.
\end{proof}

The identification of the tangent spaces of $\dbl[sm]\X$ and $\X$ allows us to extend the definition of $\opd\gamma$ to geodesics $\gamma : I \to \dbl[sm]\X$ with distinct endpoints, as long as the endpoints are identified under the canonical projection $\dbl\X \to \X$. The same formula is used to define $\opd\gamma$. Hence for any billiard loop $\beta : I \to \X$ that corresponds to a geodesic $\gamma : I \to \dbl[sm]\X$ we can define $\opd\beta = \opd\gamma$ and compute $\ker\opd\beta$ regardless of whether $\beta$ collides an odd or even number of times. We simply have to keep in mind that $\ker\opd\beta$ is related to the stability of a geodesic loop $\gamma$ only when $\beta$ has an even number of collisions.

We will construct the loops of stable figure eights using \cref{prop:Beveling}. As a consequence, their parallel defect kernels will have a very specific and convenient form as stipulated in the following lemma.

\begin{proposition}
    \label{prop:OrigamiModelParallelDefectKernelParity}
    Consider a convex polygon $\X^2 \subset \R^2$ with a billiard loop $\beta_2 : [0,1] \to \X^2$ based at $p$ such that $\beta_2'(0)$ is not parallel to $\beta_2'(1)$. Consider a sequence of $m$-polytopes $\X^m$ and billiard trajectories $\beta_m$ for $m = 3, \dotsc, n$ such that for $m \geq 2$, $\X^{m+1}$ and $\beta_{m+1}$ are derived from $\X^m$ and $\beta_m$ by applying the construction in \cref{prop:Beveling} with some choice of $k_m$ faces. Let $c_m$ be the number of collisions of $\beta_m$. Then,
    \begin{align}
        \ker\opd{\beta_n} &= U_2 \oplus V_3 \oplus \dotsb \oplus V_n, \\
        U_2 &= \vspan\{\beta_n'(0) - (-1)^{c_2} \beta_n'(1)\}, \\
        V_m &= \begin{cases}
            \vspan\{e_m\} & k_{m-1} \text{ is even,} \\
            \{0\} & k_{m-1} \text{ is odd,}
        \end{cases}
    \end{align}
    where $e_m$ denotes the $m^\text{th}$ standard basis vector of $\R^n$.
\end{proposition}
\begin{proof}
    By \cref{prop:Beveling}(1), $\beta_n(t) = (\beta_2(t),0,\dotsc,0)$ and $\beta_n'(t) = (\beta_2'(t),0,\dotsc,0)$ for $t = 0,1$. Thus let $W = \vspan\{\beta_n'(0), \beta_n'(1)\} = \R^2 \times \{(0, \dotsc, 0)\}$. By \cref{prop:Beveling}(4),
    \begin{equation}
        \label{eq:ParallelTransportDecomposition}
        P_{\beta_n} = \underbrace{P_{\beta_2}}_{P_{\beta_n}|_W} \oplus \underbrace{(-1)^{k_2} \oplus (-1)^{k_3} \oplus \dotsb \oplus (-1)^{k_{n-1}}}_{P_{\beta_n}|_{W^\perp}}.
    \end{equation}
    
    Let $\pi_t$ be as in the definition of $\opd{\beta_n}$.  \Cref{eq:ParallelTransportDecomposition} guarantees that $\vspan\{e_m\}$ is an invariant subspace of $P_{\beta_n}$, $\pi_0$, and $\pi_1$ for all $m \geq 3$. Hence,
    \begin{equation*}
        \ker\opd{\beta_n} = \ker\opd{\beta_n}|_W \oplus \ker\opd{\beta_n}|_{\vspan\{e_3\}} \oplus \dotsb \oplus \ker\opd{\beta_n}|_{\vspan\{e_n\}}.
    \end{equation*}
    \Cref{eq:ParallelTransportDecomposition} implies that $\opd{\beta_n}|_{\vspan\{e_m\}}$ is the zero map when $k_{m-1}$ is even, and is $v \mapsto -2v$ when $k_{m-1}$ is odd, which gives the formula for $V_m$.
    
    $P_{\beta_n}$ must act on $W$ by $P_{\beta_2}$, which is a reflection or rotation, depending on $\det P_{\beta_2} = (-1)^{c_2}$. The formula for $U_2$ can verified by evaluating $\opd{\beta_n}|_W$ on the orthogonal basis $(u,v)$ of $W$, where $u = \frac12(\beta_2'(0) + \beta_2'(1))$ and $v = \frac12(\beta_2'(1) - \beta_2'(0))$. The evaluation is apparent from \cref{fig:ParallelDefectKernelBeveling}, where we have written $\beta = \beta_2$ for brevity.

    \begin{figure}[h]
        \centering
        \includegraphics[width=\linewidth]{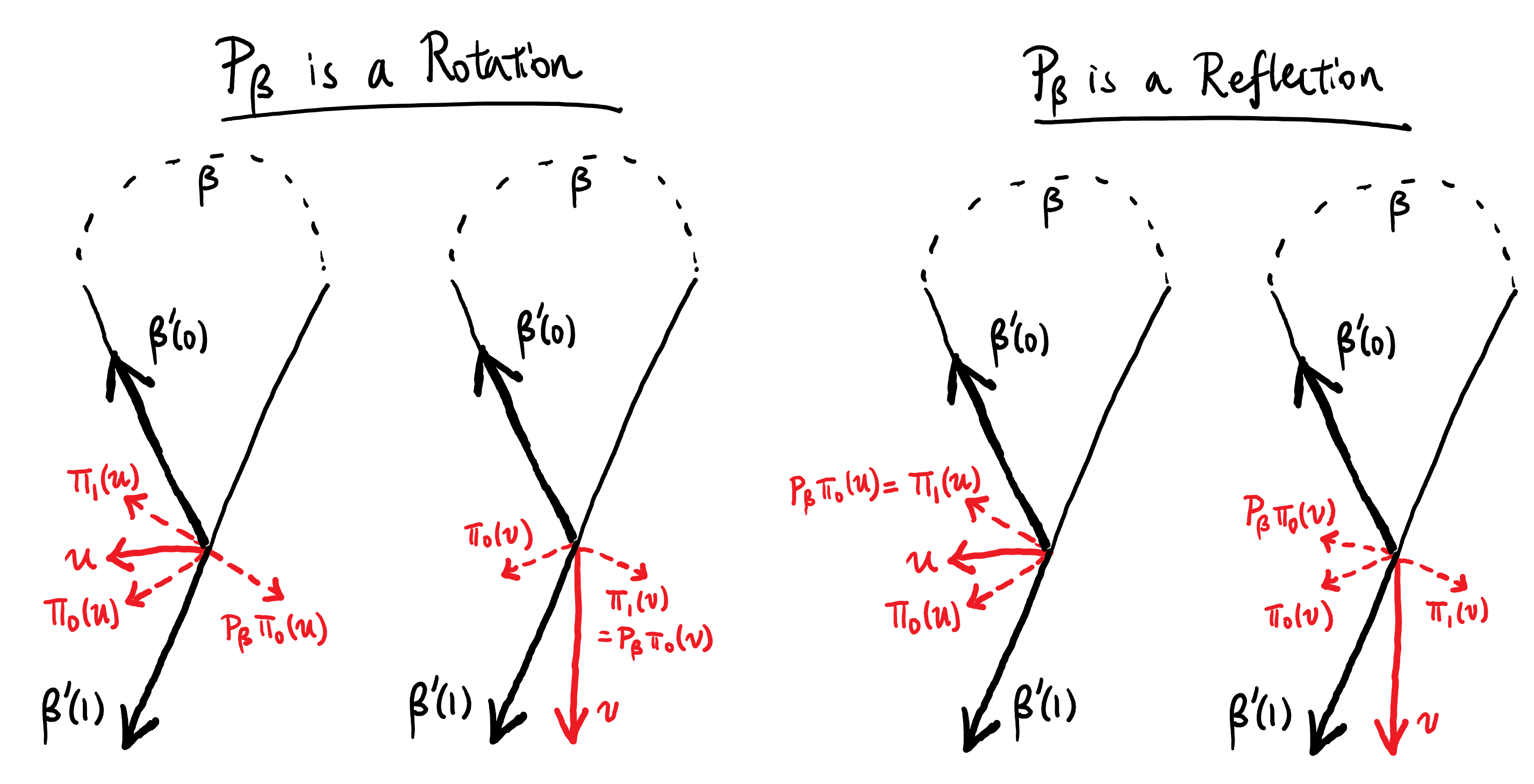}
        \caption{Computing the action of $\opd\beta$ on vectors $u$ and $v$ at the basepoint of a billiard loop $\beta$ in a convex polygon. Two cases are shown, depending on whether $P_\beta$ is a rotation or a reflection about the basepoint of $\beta$.}
        \label{fig:ParallelDefectKernelBeveling}
    \end{figure}
\end{proof}

\subsection{Product billiard trajectories}
\label{sec:ProductBilliardTraj}

In \cref{sec:Intro_Combining} we explained that part of the proof of \cref{thm:Stable2LoopPosCurv} involves combining stable figure eights in lower-dimensional manifolds into stable figure eights in higher-dimensional manifolds. The key observation is that when two billiard trajectories $\beta_i : I \to \X_i$ for $i = 1,2$ never collide simultaneously, then their \emph{product} $(\beta_1, \beta_2)$ is a billiard trajectory in $\X_1 \times \X_2$.

Given two billiard loops $\alpha, \beta : I \to \X$ that correspond to two geodesic loops in $\dbl\X$ that are based at the same point, those geodesic loops can only form a stable figure eight if they first form a \emph{stationary} geodesic bouquet: they must obey the \emph{stationarity condition} that the unit tangent vectors of the geodesics at the basepoint, pointing away from the basepoint, should sum to zero. Under the identification of the tangent spaces of $\dbl\X$ and $\X$, this translates to an analogous stationarity condition on the unit tangent vectors of $\alpha$ and $\beta$ at their basepoint. For figure eights, this condition is equivalent to requiring that $\alpha$ and $\beta$ form the same angle at their basepoint.

The following lemma gives a sufficient condition for such \emph{product billiard trajectories} to give rise to stationary figure eights.

\begin{lemma}
    \label{lem:ProductStationarity}
    For each $i = 1,2$, consider two billiard loops $\alpha_i, \beta_i : I \to \X_i$ that are based at the same point, have the same length, and satisfy the stationarity condition. Assume further that $\alpha_1$ never collides simultaneously with $\alpha_2$, and that $\beta_1$ never collides simultaneously with $\beta_2$. Then $\alpha_0 = (\alpha_1, \alpha_2)$ and $\beta_0 = (\beta_1, \beta_2)$ are billiard loops in $\X_1 \times \X_2$ that are based at the same point, have the same length, and satisfy the stationarity condition.
\end{lemma}
\begin{proof}
    For each $i = 0,1,2$, $\alpha_i$ and $\beta_i$ have the same speed. This means that the stationarity condition for $\alpha_i$ and $\beta_i$ is equivalent to the constraint that $\alpha_i'(0) - \alpha_i'(1) + \beta_i'(0) - \beta_i'(1) = 0$. The lemma then follows by basic geometry.
\end{proof}

The following lemma will be useful to relate the stability of the figure eights that arise from product billiard trajectories to the stability of figure eights that come from the trajectories in each factor.

\begin{lemma}
    \label{lem:ParallelDefectKernelProduct}
    For $i = 1,2$, let $\beta_i : I \to \X_i$ be a billiard loop in a $n_i$-polytope $\X_i$. If $\beta_1$ and $\beta_2$ never collide simultaneously, then
    \begin{equation}
        \ker \opd{(\beta_1, \beta_2)} \subset \ker\opd{\beta_1} \times \ker\opd{\beta_2}.
    \end{equation}
\end{lemma}
\begin{proof}
    For a billiard trajectory $\beta : I \to \X$, define $\pi_t^\beta : T_{\beta(t)}\X \to T_{\beta(t)}\X$ to be the projection onto $\beta'(t)^\perp$.
    
    Suppose that $u = (u_1,u_2) \in \ker \opd{\beta}$, where $\beta = (\beta_1,\beta_2)$ and $u_i \in T_{\beta_i(0)}\X_i$. Observe that $\beta'(t)^\perp = \beta_1'(t)^\perp \times \beta_2'(t)^\perp \oplus \vspan\{r_t\}$, where $r_t = (-s_2^2\beta_1'(t), s_1^2\beta_2'(t))$ and $s_i$ is the speed of $\beta_i$. Thus for $t = 0,1$ we have the direct sum decomposition
    \begin{equation}
        \label{eq:ParallelDefectKernelProduct_ExpandOrthoProj}
        \pi_t^\beta u = (\pi_t^{\beta_1} u_1, \pi_t^{\beta_2} u_2) + \lambda_t r_t \qquad \text{for some } \lambda_t \in \R.
    \end{equation}
    It can be verified that $P_\beta = (P_{\beta_1}, P_{\beta_2})$. Since $u \in \ker\opd\beta$, we have $P_\beta\pi_0^\beta u = \pi_1^\beta u$. Expanding both sides using \cref{eq:ParallelDefectKernelProduct_ExpandOrthoProj} and the fact that $P_\beta r_0 = r_1$ yields $(P_{\beta_1} \pi_0^{\beta_1} u_1, P_{\beta_2} \pi_0^{\beta_2} u_2) + \lambda_0 r_1 = (\pi_1^{\beta_1} u_1, \pi_1^{\beta_2} u_2) + \lambda_1 r_1$. The uniqueness of the direct sum decomposition implies that $P_{\beta_i} \pi_0^{\beta_i} u_i = \pi_1^{\beta_i} u_i$ for $i = 1,2$. In other words, $u \in \ker\opd{\beta_1} \times \ker\opd{\beta_2}$.
\end{proof}

\subsection{The construction of stable figure eights}
\label{sec:Stable2Loop}

To prove \cref{thm:Stable2LoopPosCurv}, we begin by constructing a 3-polytope, 4-polytope, and 5-polytope whose doubles contain stable figure eights. We will frequently apply the construction from \cref{prop:Beveling} to eventually get a billiard loop $\tilde\beta : I \to \tilde\X$ that corresponds to a geodesic loop in the final stable figure eight.

Recall that by \cref{eq:Beveling}, $\tilde\X$ is the intersection of $\X \times \R$ with wedges $\W_{h_i}(\HS_i)$, where $\HS_1, \dotsc, \HS_k$ are halfspaces bounding $\X$ whose boundaries are encountered in sequence by $\beta$. The proof of \cref{prop:Beveling} allows to choose particular values of the $h_i$'s, as long as they are sufficiently small in magnitude and their alternating sum vanishes: $h_1 - h_2 + \dotsb + (-1)^{k-1}h_k = 0$. (To simplify notation, if $F$ is a face of an $n$-polytope $\X$ that lies in the boundary of a half-space $\HS$ containing $\X$, then let $\W_h(F)$ stand for $\W_h(\HS)$.)

\Cref{prop:OrigamiModelParallelDefectKernelParity} shows that if we wish to iterate the construction from \cref{prop:Beveling} to produce geodesic loops with parallel defect kernels of low dimension, then it is important to carefully choose the parities of several integers. For instance, the parity of $c_2$ controls $U_2$. It would also be ideal for each application of \cref{prop:Beveling} to involve an odd number of faces so that as many of the $V_m$'s are 0-dimensional as possible; but we also need $c_n = c_2 + k_2 + \dotsb + k_{n-1}$ to be even so that the resulting billiard loop $\beta_n$ will correspond to a geodesic loop in $\dbl\X^n$. Hence we should pay attention to the \emph{parity sequence} of a series of applications of \cref{prop:Beveling}, by which we mean the sequence of parities of $c_2, k_2, \dotsc, k_{n-1}$.

Let us proceed to perform the main construction for the 3-dimensional case of \cref{thm:Stable2LoopPosCurv}.

\begin{figure}[p]
    \centering
    \includegraphics[width=0.8\textwidth]{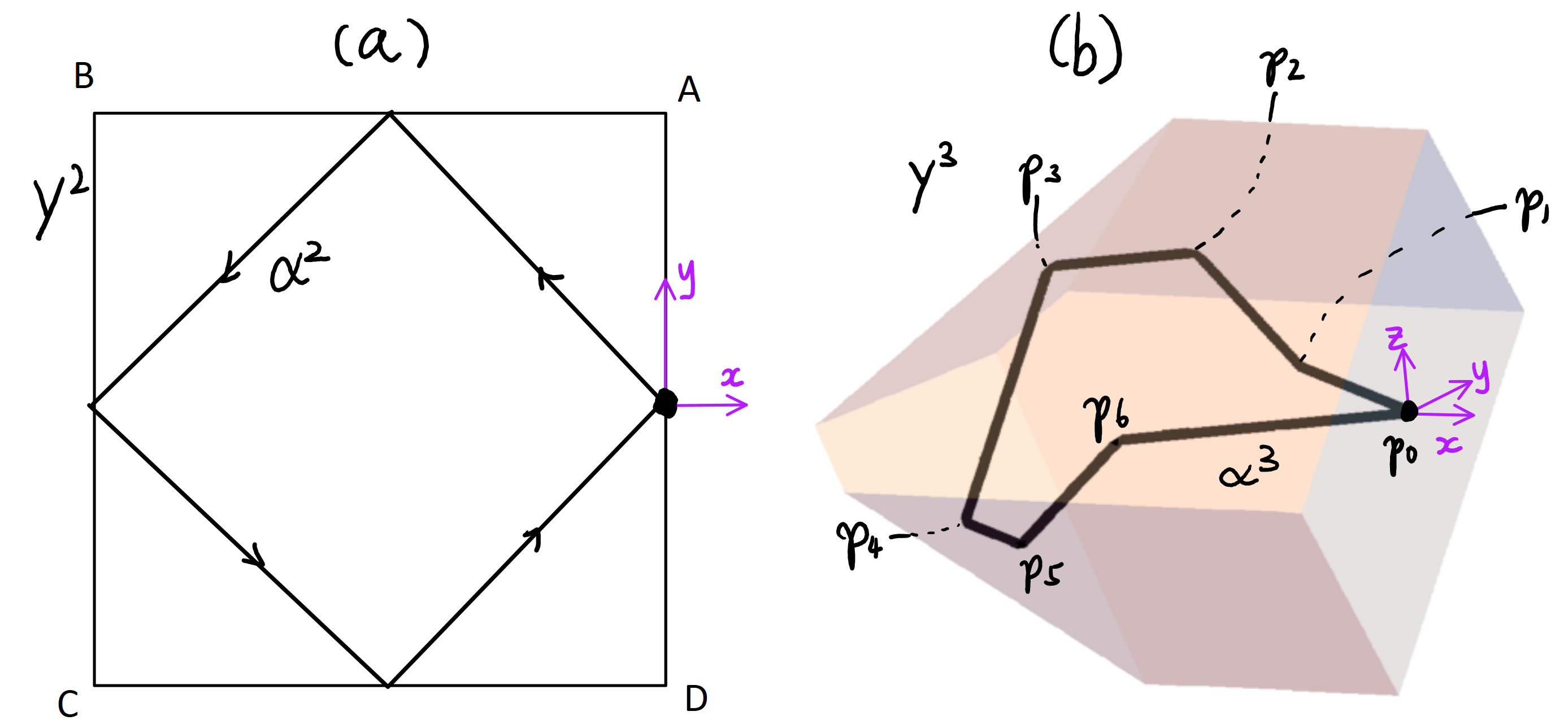}
    \includegraphics[width=0.8\textwidth]{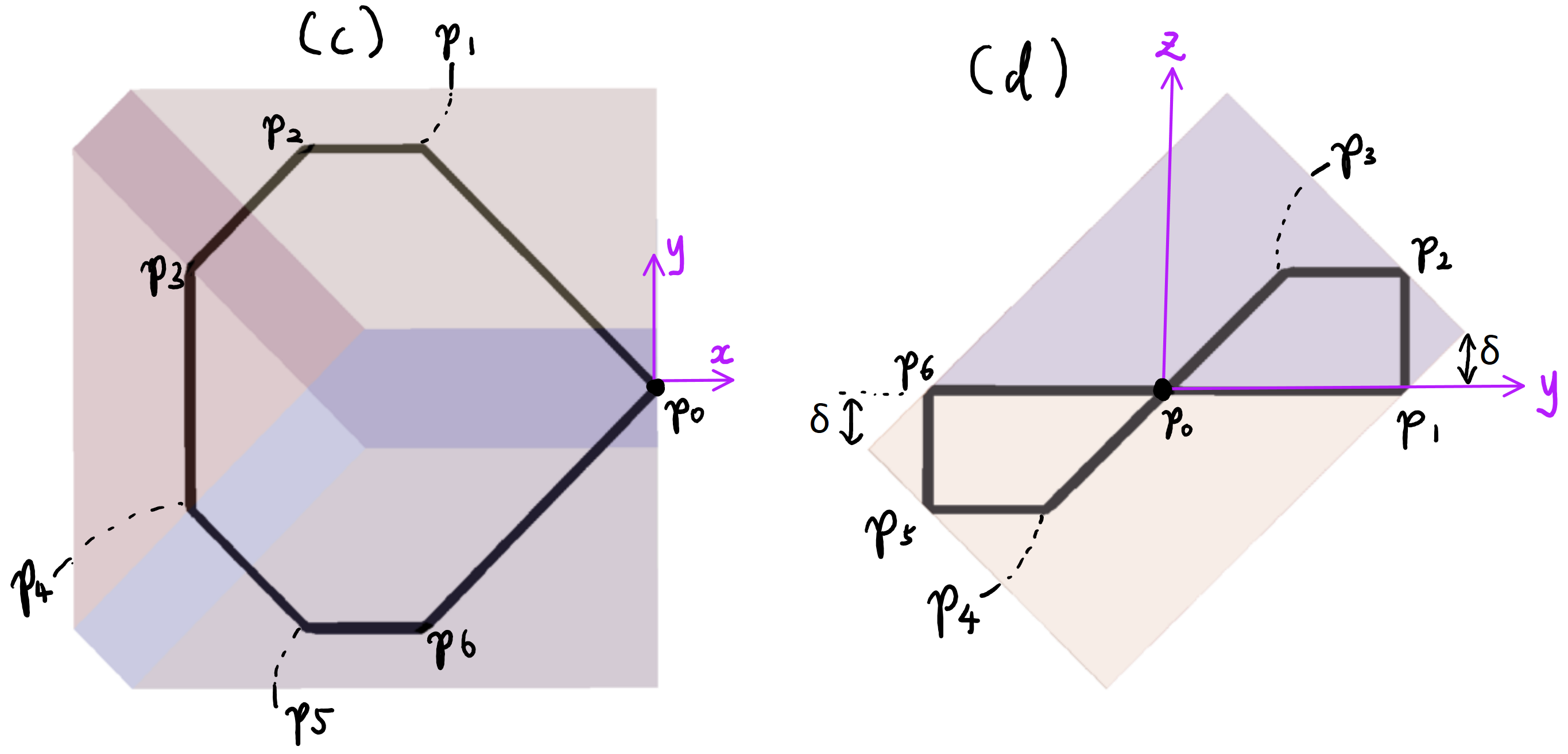}
    \includegraphics[width=0.8\textwidth]{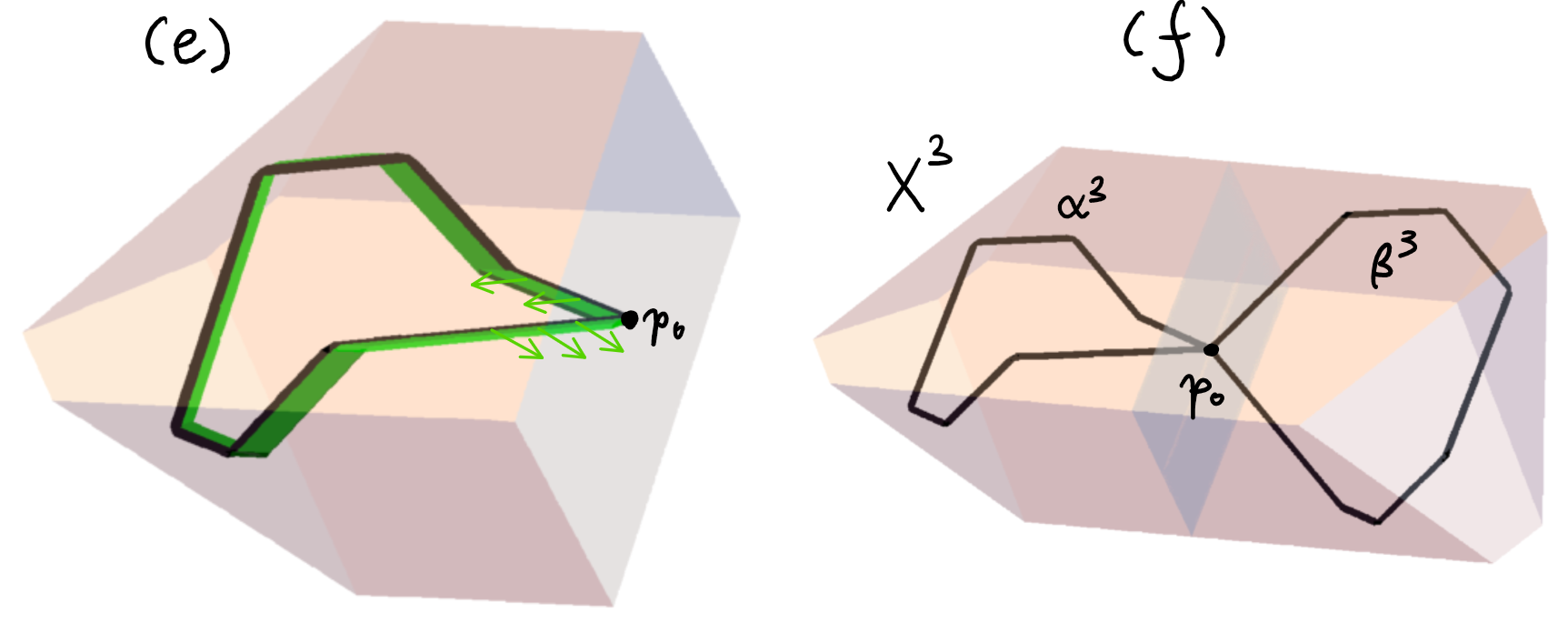}
    \caption{Various stages from the construction of a stable figure eight in the proof of \cref{prop:StableFigureEightFlat3D}. (a) The initial polygon $\Y^2$ and billiard trajectory $\alpha_2$. (b)--(d) Three views of $\Y^3$ and $\alpha_3$. (e) The green strip represents a projection of a parallel vector field along the geodesic loop corresponding to $\alpha_3$. A few arrows from that field are drawn. (f) $\X^3$ and the billiard trajectories that correspond to the stable figure eight in \cref{prop:StableFigureEightFlat3D}.}
    \label{fig:StableFigureEightFlat3D}
\end{figure}

\begin{proposition}
    \label{prop:StableFigureEightFlat3D}
    There exists a 3-polytope $\X^3 \subset \R^3$ and simple billiard loops $\alpha_3$ and $\beta_3$ in $\X^3$ of equal length that correspond to two geodesic loops of a simple and irreducible stable figure eight $G^3$ in $\dbl\X^3$. Moreover, $\alpha_3$ only intersects $\beta_3$ at their common basepoint.
\end{proposition}
\begin{proof}
    We want to find such $\alpha_3, \beta_3$ that correspond to geodesic loops $\gamma_1, \gamma_2$. To obtain each geodesic loop, we will apply the \cref{prop:Beveling} once using the parity sequence ``odd, odd.''

    Now consider the square $\Y^2  = [-2,0] \times [-1,1] \subset \R^2$, whose vertices are $A$, $B$, $C$, and $D$, listed counterclockwise from $A = (0,1)$. Let $\alpha_2 : [0,1] \to \Y^2$ be the billiard trajectory that starts from the origin, collides at the points $(-1,1)$, $(-2,0)$ and $(-1,-1)$ in order, and then returns to the origin (see \cref{fig:StableFigureEightFlat3D}(a)). Let us apply \cref{prop:Beveling} to get
    \begin{equation*}
        \Y^3 = \Y^2 \times \R \cap \W_{\delta}(AB) \cap \W_0(BC) \cap \W_{-\delta}(CD),
    \end{equation*}
    for sufficiently small $\delta > 0$. This gives a billiard trajectory $\alpha_3$ in $\Y^3$, both of which are depicted from three different points of view in \cref{fig:StableFigureEightFlat3D}(b)--(d). \Cref{prop:Beveling}(1) guarantees that $\alpha_3(1) = \alpha_3(0) = (0,0,0)$ and $\alpha_3'(t) = (\alpha_2'(t), 0)$ for $t = 0,1$.

    Now let $\rho \in SO(3)$ be the rotation by angle $\pi$ about the line spanned by the vector $(0,1,1)$. One can verify that $\X^3 = \Y^3 \cup \rho(\Y^3)$ is a 3-polytope. As $\alpha_3$ is a billiard loop in $\X^3$ with an even number of collisions, it corresponds to a geodesic loop $\gamma_1$ in $\dbl\X^3$. Similarly, $\beta_3 = \rho \circ \alpha_3$ also corresponds to a geodesic loop $\gamma_2$ in $\dbl\X^3$ that shares the same basepoint $p$ as $\gamma_1$. $\X^3$, $\alpha_3$ and $\beta_3$ are illustrated in \cref{fig:StableFigureEightFlat3D}(f).
    
    It can be verified that $\gamma_1$ and $\gamma_2$ form an irreducible stationary geodesic bouquet $G^3$ in $\dbl\X^3$ that is based at $p$. (The stationarity condition is satisfied as $\alpha_3$ and $\beta_3$ form the same angle at their basepoint.) \Cref{prop:OrigamiModelParallelDefectKernelParity} implies that $\ker\opd{\gamma_1} = \ker\opd{\alpha_3}$ is the $y$-axis.\footnote{This can be independently checked by observing that the parallel transport map is rotation by angle $\pi$ around the $y$-axis. This is consistent with the illustration in \cref{fig:StableFigureEightFlat3D}(e) of the projection of the parallel vector field along $\gamma_1$ that is induced by the quotient map $\dbl\X^3 \to \X^3$.} We also have that $\ker\opd{\gamma_2} = \rho(\ker\opd{\alpha_3})$ is the $z$-axis, which intersects $\ker\opd{\gamma_1}$ only at the origin. Therefore $G^3$ is stable by \cref{cor:NullvarsToIntersectKernels_Flat}.
    
    \Cref{prop:Beveling}(5) implies that $\alpha_3$ and $\beta_3$ are simple. By construction, they intersect only at their common basepoint. As a result, $G^3$ is also simple.
\end{proof}

\begin{remark}
    Let us give motivation for this construction. \label{rem:StableFigureEightFlat3D}
    In the proof of \cref{prop:StableFigureEightFlat3D}, the rotation $\rho$ performs two roles. It rotates $\alpha_3$ to $\beta_3$ while also rotating $\ker\opd{\gamma_1}$ to $\ker\opd{\gamma_2}$, with the result that $\ker\opd{\gamma_1} \cap \ker\opd{\gamma_2} = \{0\}$. The rotation also ensures that the tangent vectors of the resulting stable figure eight do not all lie in the same plane, so the stable figure eight is not a closed geodesic. The heights of the wedges are chosen so that the face $F$ of $\Y^3$ that lies on the $yz$-plane is the rectangle shown in \cref{fig:StableFigureEightFlat3D}(d), which is mapped back to itself by $\rho$. This helps $\Y^3 \cup \rho(\Y^3)$ to be a convex polyhedron. This fact is also helped by the choice of $\Y^2$ as a square, so that the dihedral angles in $\Y^3$ between $F$ and its adjacent faces would be $\pi/2$.
\end{remark}

Next, let us perform the main construction for the 4-dimensional case of \cref{thm:Stable2LoopPosCurv}. 

\begin{proposition}
    \label{prop:StableFigureEightFlat4D}
    There exists a 4-polytope $\X^4$ and simple billiard loops $\alpha_4$ and $\beta_4$ in $\X^4$ of equal length that correspond to two geodesic loops of a simple and irreducible stable figure eight $G^4$ in $\dbl\X^4$. Moreover, $\alpha_4$ only intersects $\beta_4$ at their common basepoint.
\end{proposition}
\begin{proof}
    The desired geodesic loops will correspond to billiard trajectories $\alpha_4$ and $\beta_4$ in a convex 4-polytope $\X^4$. $\alpha_4$ and $\beta_4$ are each obtained from two applications of \cref{prop:Beveling}, using the respective parity sequences ``odd, even, odd'' and ``even, odd, odd''. \comment{We will motivate these choices in a remark after the proof.}
    
    We begin with convex polygons $\Y^2, \ZZ^2 \subset \R^2$ that contain billiard trajectories $\alpha_2$ and $\beta_2$ respectively, as illustrated in \cref{fig:StableFigureEightFlat4D}(a) and (d). The $x$-coordinates of the vertices $B$ and $F$ in \cref{fig:StableFigureEightFlat4D}(a) depend on a parameter $1 < \lambda < 3$, but once $\lambda$ has been fixed, it uniquely determines the positions of vertices $C$ and $E$ such that $\alpha_2$ is a billiard trajectory in $\Y^2$ with the angles labeled as shown. Similarly, a parameter $1 < \mu < 3$ determines the positions of vertices $H$ and $J$ in \cref{fig:StableFigureEightFlat4D}(d), which then uniquely determines the position of the vertex $I$.
    
    Let us choose $\lambda$ and $\mu$ so that $\alpha_2$ and $\beta_2$ will have equal length. Referring to points labeled in \cref{fig:StableFigureEightFlat4D}(a) and (d), the equal-length criterion is equivalent by some trigonometry to:
    \begin{equation}
        \label{eq:StableFigureEightFlat4D_EqualLength}
        \underbrace{2(1 + \sin\tfrac{\pi}4/\sin\tfrac{3\pi}8)}_{\approx 3.53}P_1P_2 = \underbrace{(2 + \sqrt2)}_{\approx 3.41}Q_1Q_2.
    \end{equation}
    As $\lambda$ ranges from 1 to 3, the distance $P_1P_2$ ranges from $\sqrt2\tan\frac{3\pi}{16}$ to $\sqrt2$. Similarly, as $\mu$ ranges from 1 to 3, $Q_1Q_2$ ranges from $\sqrt2\tan\frac\pi8$ to $\sqrt2$. Thus we could choose $\lambda$ and $\mu$ so that $Q_1Q_2 \approx \sqrt2$ and $P_1P_1 \in (\sqrt2\tan\tfrac{3\pi}{16}, Q_1Q_2)$ to satisfy \cref{eq:StableFigureEightFlat4D_EqualLength}.

    In accordance with the respective parity sequences, $\alpha_2$ collides an odd number of times, and $\beta_2$ collides an even number of times.
    
    \begin{figure}[p]
        \centering
        \includegraphics[width=0.9\textwidth]{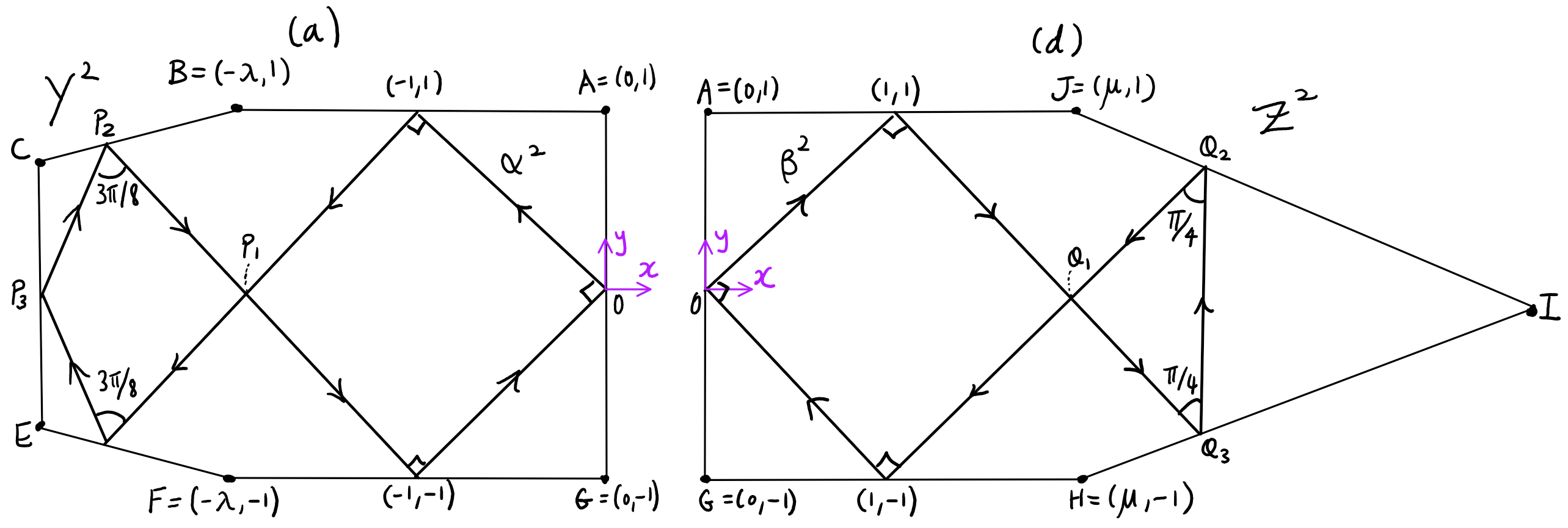}
        \includegraphics[width=0.85\textwidth]{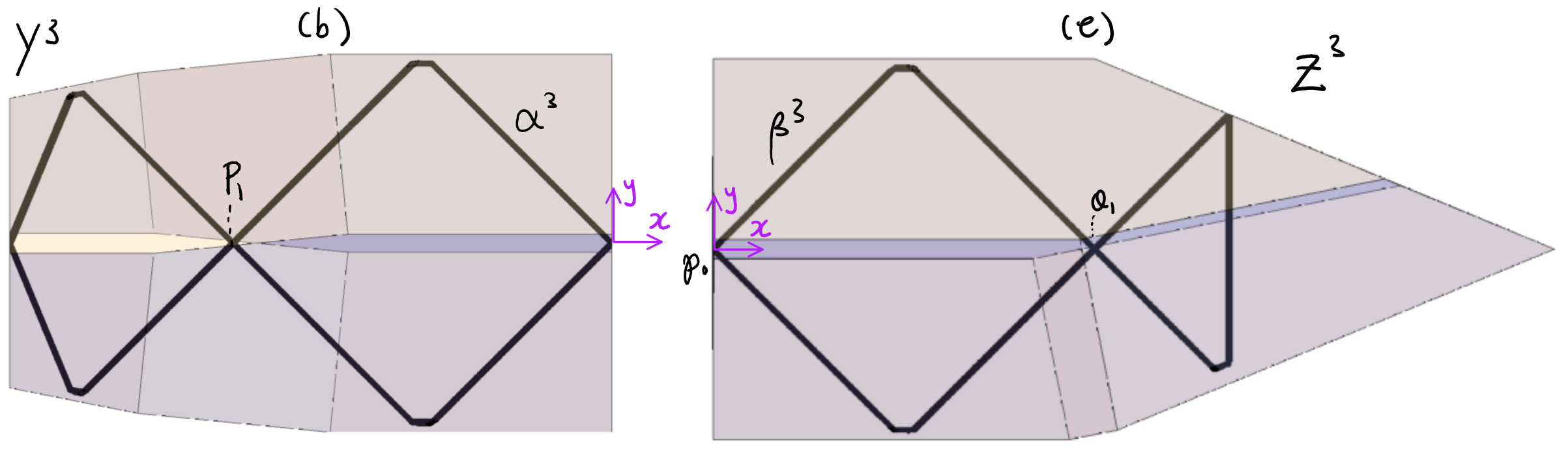}
        \includegraphics[width=0.9\textwidth]{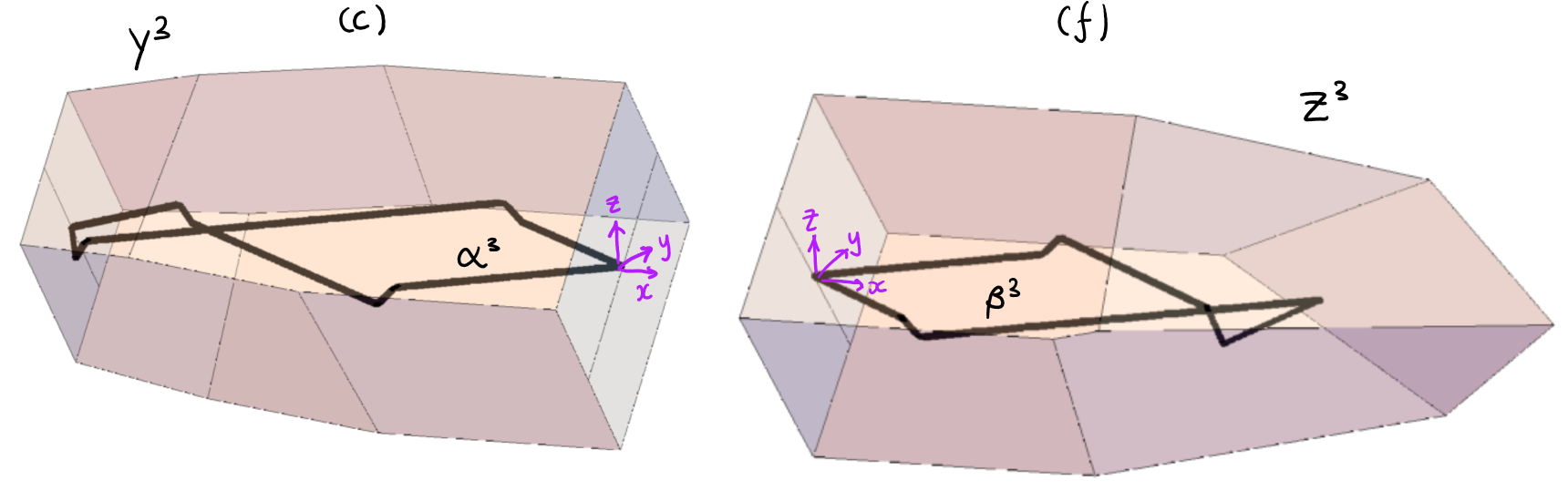}
        \includegraphics[width=0.7\textwidth]{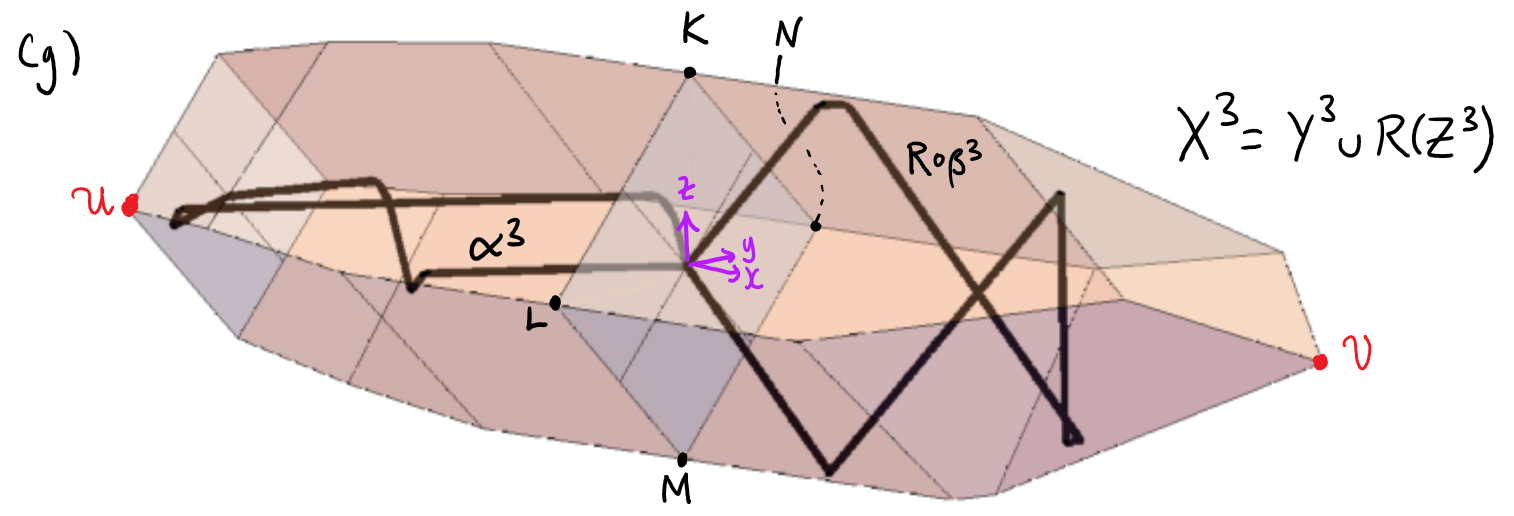}
        \caption{Various stages from the construction of a stable figure eight in the proof of \cref{prop:StableFigureEightFlat4D}. (a) The polygon $\Y^2$ and billiard trajectory $\alpha_2$. (b)--(c) Two views of $\Y^3$ and $\alpha_3$. (d) The polygon $\ZZ^2$ and billiard trajectory $\beta_2$. (e)--(f) Two views of $\ZZ^3$ and $\beta_3$. (g) $\X^3$ and the billiard trajectories $\alpha_3$ and $R \circ \beta_3$ that correspond to the stable figure eight in \cref{prop:StableFigureEightFlat4D}.}
        \label{fig:StableFigureEightFlat4D}
    \end{figure}
    
    Let $\delta$ be a small positive number. Apply \cref{prop:Beveling} to $\Y^2$ and $\alpha_2$ to get a 3-polytope,
    \begin{equation*}
        \Y^3 = \Y^2 \times \R \cap \W_\delta(AB) \cap \W_\delta(EF) \cap \W_{-\delta}(BC) \cap \W_{-\delta}(FG),
    \end{equation*}
    that contains a billiard trajectory $\alpha_3$. They are illustrated in \Cref{fig:StableFigureEightFlat4D}(b)--(c). Next, apply \cref{prop:Beveling} to $\ZZ^2$ and $\beta_2$ to get a 3-polytope,
    \begin{equation*}
        \ZZ^3 = \ZZ^2 \times \R \cap \W_\delta(JA) \cap \W_0(HI) \cap \W_{-\delta}(GH),
    \end{equation*}
    that contains a billiard trajectory $\beta_3$ (``odd'' in the parity sequence). \Cref{fig:StableFigureEightFlat4D}(e)--(f) depicts $\ZZ^3$ and $\beta_3$. As a consequence of the theorem, $\alpha_3$ collides at most once with each face of $\Y^3$, and similarly for $\beta_3$ in $\ZZ^3$.
    
    Note that $\alpha_3$ is simple: its projection onto the $xy$-plane intersects itself at the point $P_1$ as shown in \cref{fig:StableFigureEightFlat4D}(b), but the two segments of $\alpha_3$ that are involved lie in the planes $\{z = 2\delta\}$ and $\{z = -2\delta\}$ respectively, so they cannot intersect. Similarly, $\beta_3$ is simple: its projection onto the $xy$-plane intersects itself at the point $Q_1$ as shown in \cref{fig:StableFigureEightFlat4D}(e), but the two segments involved lie in the planes $\{z = 2\delta\}$ and $\{z = -2\delta\}$ respectively, so they cannot intersect.
    
    The same value of $\delta$ is used to construct $\Y^3$ and $\ZZ^3$ so that the rectangle $KLMN$ (see \cref{fig:StableFigureEightFlat4D}(g)) is a face of both $\Y^3$ and $\ZZ^3$. Now let $R \in O(3)$ be the reflection that fixes the $x$ coordinate but interchanges the $y$ and $z$ coordinates. Observe that $R$ maps $KLMN$ back to itself, so that $\X^3 = \Y^3 \cup R(\ZZ^3)$ is a convex polyhedron, and that $R \circ \beta_3$ is a billiard trajectory in $\X^3$. \Cref{fig:StableFigureEightFlat4D}(g) illustrates $\X^3$, $\alpha_3$ and $R \circ \beta_3$. This transformation $R$ is applied to ensure that the the tangent vectors of the final billiard trajectories $\alpha_4$ and $\beta_4$ at the origin will not lie in the same plane, so that the corresponding geodesic loops do not form a closed geodesic. Eventually this will lead to the irreducibility of the final stable figure eight.\footnote{The rotation performed in the proof of \cref{prop:StableFigureEightFlat3D} also fulfills this purpose.}
    
    Next, we will construct $\alpha_4$ and $\beta_4$ from $\alpha_3$ and $R \circ \beta_3$ by applying \cref{prop:Beveling} again. Three of the faces of $\Y^3$ that $\alpha_3$ collides with are portions of $CE \times \R$ and $\partial \W_\delta(EF)$; in \cref{fig:StableFigureEightFlat4D}(g), they are the three faces around the vertex $u$. Index those faces $A_1$, $A_2$ and $A_3$ in the order that $\alpha_3$ collides with them. Similarly, three of the faces that $R \circ \beta_3$ collides with are images under $R$ of the faces of $\ZZ^3$ that are portions of $\partial\W_0(HI)$ and $IJ \times \R$; in \cref{fig:StableFigureEightFlat4D}(g), they are the three faces around the vertex $v$. Index those faces $B_1$, $B_2$ and $B_3$ in the order that $\beta_3$ collides with them. Choose some sufficiently small $\varepsilon > 0$ and apply \cref{prop:Beveling} to $\X^3$ and the billiard trajectories $\alpha_3$ and $\beta_3$ get 4-polytopes,
    \begin{align*}
        \X^4_\alpha &= \X^3 \times \R
        \cap \W_\varepsilon(A_1) \cap \W_0(A_2) \cap \W_{-\varepsilon}(A_3)
        \\
        \X^4_\beta &= \X^3 \times \R \cap \W_\varepsilon(B_1) \cap \W_0(B_2) \cap \W_{-\varepsilon}(B_3),
    \end{align*}
    that contain billiard trajectories $\alpha_4$ and $\beta_4$ respectively. (This gives the final ``odd'' in each parity sequence.) Finally, observe that $\alpha_4$ and $\beta_4$ are also billiard trajectories in $\X^4 = \X^4_\alpha \cap \X^4_\beta$.
    
    \Cref{prop:Beveling}(2) implies that $\alpha_4$ and $\beta_4$ each collide an even number of times, so they correspond to geodesic loops $\gamma_1$ and $\gamma_2$ in $\dbl\X^4$ that are based at the same point $p$, giving a stationary figure eight $G^4$. (To verify the stationarity condition, observe that $\alpha_2$ and $\beta_2$ form the same angle at their basepoint. This property is preserved in $\alpha_3$ and $\beta_3$ and also in $\alpha_4$ and $\beta_4$ due to \cref{prop:Beveling}(1).) Given the parity sequence used to construct $\alpha_4$, \cref{prop:OrigamiModelParallelDefectKernelParity} implies that $\ker\opd{\gamma_1} = \ker\opd{\alpha_4} = \{0\} \times \R^2 \times \{0\}$. A similar application of the arguments used to prove \cref{prop:OrigamiModelParallelDefectKernelParity} to the parity sequence used to construct $\beta_4$, but taking into account the reflection $R$, implies that $\ker\opd{\gamma_2} = \R \times \{(0,0,0)\}$. Therefore $\ker\opd{\gamma_1} \cap \ker\opd{\gamma_2} = \{0\}$, and $G^4$ is stable by \cref{cor:NullvarsToIntersectKernels_Flat}.
    
    $\alpha_4$ and $\beta_4$ are simple by \cref{prop:Beveling}(5). They also only intersect at their common basepoint, because $\im({\pr_3} \circ \alpha_4)$ lies in the convex hull of $\im(\alpha_3)$, and similarly for $\beta_4$. Thus $G^4$ is simple. $\alpha_2$ and $\beta_2$ have the same length, so \cref{prop:Beveling}(\ref{enum:Beveling_Length}) implies that $\alpha_i$ has the same length as $\beta_i$ for $i = 3,4$ as well.
\end{proof}

The main construction in the proof of the 5-dimensional case of \cref{thm:Stable2LoopPosCurv} will use product billiard trajectories.

\begin{proposition}
    \label{prop:StableFigureEightFlat5D}
    There exists a 5-polytope $\X^5$ and simple billiard loops $\alpha_5$ and $\beta_5$ in $\X^5$ of equal length that correspond to two geodesic loops of a simple, irreducible and stable figure eight $G^5$ in $\dbl\X^5$. Moreover, $\alpha_5$ and $\beta_5$ intersect only at their common basepoint.
\end{proposition}
\begin{proof}
    We will begin with the 3-polytope $\X^3$ and the billiard trajectories $\alpha_3$ and $R \circ \beta_3$ in it, from the proof of \cref{prop:StableFigureEightFlat4D}. $\alpha_3$ and $\beta_3$ have equal length, and were derived using \cref{prop:Beveling} from other billiard trajectories $\alpha_2 : I \to \Y^2$ and $\beta_2 : I \to \ZZ^2$. $\beta_3$ is composed with a reflection $R$ that interchanges the $y$ and $z$ coordinates to get $R \circ \beta_3$. We will also consider a regular decagon $\X^2$ and two periodic billiard trajectories $\phi, \psi : I \to \X^2$ that form regular pentagons in it, starting and ending at the same point $y$ (see \cref{fig:StableFigureEightFlat5D}), and have equal length. All of the above billiard trajectories have an odd number of collisions.
    
    \begin{figure}[h]
        \centering
        \includegraphics[width=5cm]{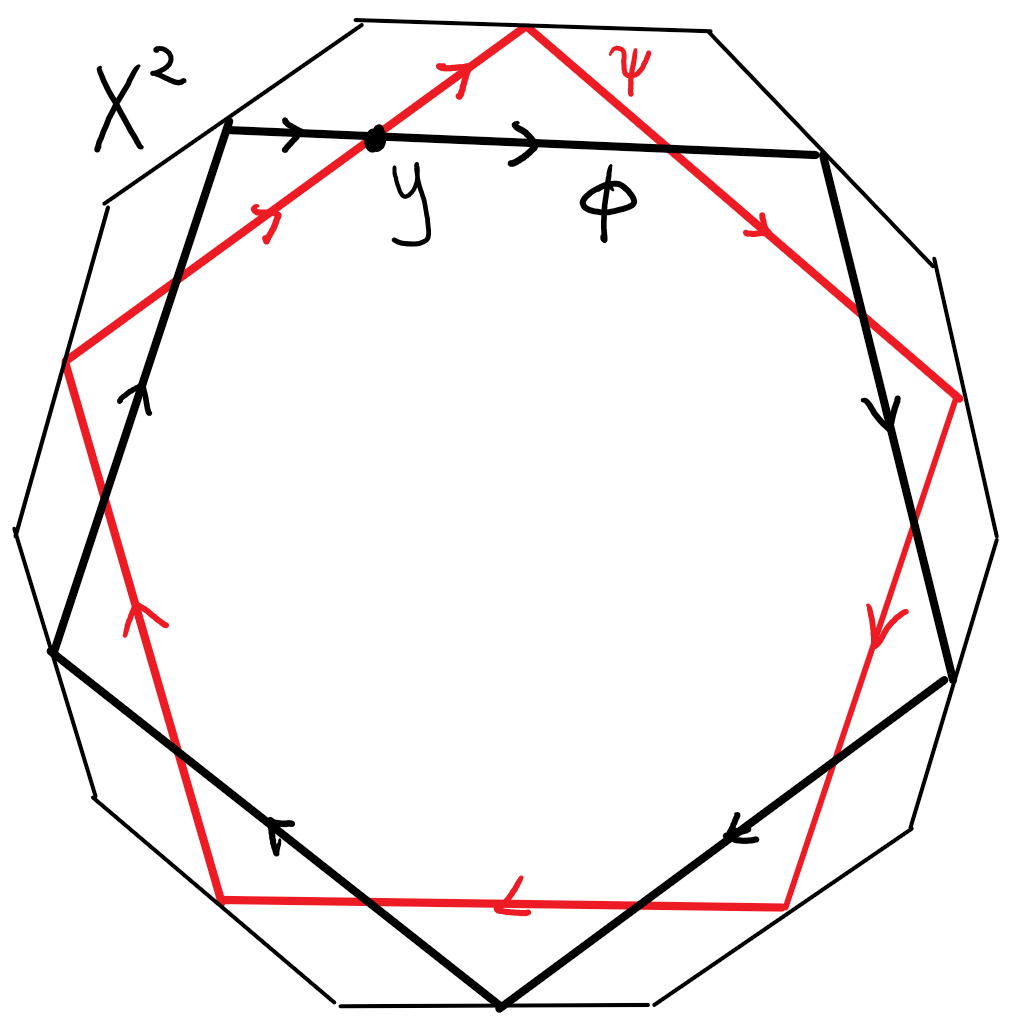}
        \caption{A regular decagon $\X^2$ with two periodic billiard trajectories, $\phi$ (black) and $\psi$ (red), starting and ending at a point $y \in \X^2$.}
        \label{fig:StableFigureEightFlat5D}
    \end{figure}
    
    \Cref{prop:Beveling}(6) allows us to assume that the collision times of $\alpha_3$ will be arbitrarily close to those of $\alpha_2$. Hence, $\alpha_3$ will never collide simultaneously with $\phi$, and so $\alpha_5 = (\alpha_3, \phi)$ is a billiard trajectory in $\X^3 \times \X^2$ (see the discussion in \cref{sec:Intro_Combining}). $\alpha_5$ also has an even number of collisions, so it corresponds to a geodesic loop $\gamma_1$ in $\dbl(\X^3 \times \X^2)$. Similarly, we may assume that $\beta_5 = (R \circ \beta_3, \psi)$ is a billiard trajectory in $\X^3 \times \X^2$ with an even number of collisions, which corresponds to a geodesic loop $\gamma_2$ in $\dbl(\X^3 \times \X^2)$. $\gamma_1$ and $\gamma_2$ form a stationary figure eight $G^5$, because \cref{lem:ProductStationarity} implies that the stationarity condition is satisfied. Another implication is that $\alpha_5$ has the same length as $\beta_5$.

    Next, we study the stability of $G^5$  by computing parallel defect kernels. From \cref{prop:OrigamiModelParallelDefectKernelParity} we get $\ker\opd{\alpha_3} = \{0\} \times \R^2$ and $\ker\opd{R \circ \beta_3} = R(\ker\opd{\beta_3}) = R(\R \times \{(0,0)\}) = \R \times \{(0,0)\}$. $P_\phi$ is the product of an odd number of reflections, so it has to be the reflection that fixes the tangent vector at the basepoint. Thus we can apply the definition of parallel defect kernels to show that $\ker\opd\phi = \vspan\{\phi'(0)\}$. Similarly, $\ker\opd\psi = \vspan\{\psi'(0)\}$. By \cref{lem:ParallelDefectKernelProduct},
    \begin{multline*}
        \ker \opd{\alpha_5} \cap \ker \opd{\beta_5} \subset (\ker \opd{\alpha_3} \times \ker\opd\phi) \cap (\ker \opd{R \circ \beta_3} \times \ker\opd\psi)
        \\
        = (\ker \opd{\alpha_3} \cap \ker \opd{R \circ \beta_3}) \times (\ker\opd\phi \cap \ker\opd\psi) = \{0\}.
    \end{multline*}
    Therefore $G^5$ is stable. Its irreducibility can also be verified. 

    $\alpha_3$ and $\beta_3$ are simple as shown in the proof of \cref{prop:StableFigureEightFlat4D}. Moreover, $\alpha_3$ and $R \circ \beta$ intersect only at their basepoints. Thus $\alpha_5$ and $\beta_5$ are simple billiard loops that intersect only at their basepoints. As a result, $G^5$ is also simple.
\end{proof}

Finally we can combine all of our previous results into a proof of our main result, \cref{thm:Stable2LoopPosCurv}.

\begin{proof}[Proof of \cref{thm:Stable2LoopPosCurv}]
    From \cref{prop:StableFigureEightFlat3D,prop:StableFigureEightFlat4D,prop:StableFigureEightFlat5D} we obtained the convex polytopes $\X^3$, $\X^4$, and $\X^5$ in dimensions 3, 4 and 5 whose doubles contain the stable figure eights $G^3$, $G^4$ and $G^5$ respectively. $G^n$ was derived from a pair of simple billiard loops $\alpha_n$ and $\beta_n$ of equal length. Now fix some $n \in \{3,4,5\}$ and integer $m \geq 1$. Since $\alpha_3$ and $\beta_3$ were obtained by applying \cref{prop:Beveling}, property (6) of that proposition allows us to use different choices of parameters in the proof of the proposition to produce $m$ new versions of $\alpha_3$ and $\beta_3$ in $m$ convex polyhedra similar to $\X^3$, denoted by $\alpha_3^{(i)}, \beta_3^{(i)} : I \to \X^3_{(i)}$ for $i = 1, \dotsc, m$, so that no two of the billiard trajectories $\alpha_n, \alpha_3^{(1)}, \dotsc, \alpha_3^{(m)}$ ever collide simultaneously. Then $\alpha = (\alpha_n, \alpha_3^{(1)}, \dotsc, \alpha_3^{(m)})$ is a billiard trajectory in $\X^{n + 3m} = \X^n \times \X^3_{(1)} \times \dotsb \X^3_{(m)}$. Similarly, $\beta = (\beta_n, \beta_3^{(1)}, \dotsc, \beta_3^{(m)})$ is a billiard trajectory in $\X^{n + 3m}$.
    
    $\alpha$ and $\beta$ obey the stationarity condition by \cref{lem:ProductStationarity}, and correspond to geodesics $\gamma_1$ and $\gamma_2$ that form a stationary figure eight $G^{n + 3m}$ in $\dbl\X^{n + 3m}$. From the proofs of \cref{prop:StableFigureEightFlat3D,prop:StableFigureEightFlat4D,prop:StableFigureEightFlat5D} we know that $\ker\opd{\alpha_n} \cap \ker\opd{\beta_n} = \{0\}$ and $\ker\opd{\alpha_3^{(i)}} \cap \ker\opd{\beta_3^{(i)}} = \{0\}$. Then \cref{lem:ParallelDefectKernelProduct} implies that
    \begin{multline*}
        \ker\opd{\gamma_1} \cap \ker\opd{\gamma_2} = \ker\opd{\alpha} \cap \ker\opd{\beta}
        \\
        \subset \left( \ker\opd{\alpha_n} \times \prod_{i = 1}^m \ker\opd{\alpha_3^{(i)}} \right) \cap \left( \ker\opd{\beta} \times \prod_{i = 1}^m \ker\opd{\beta_3^{(i)}} \right) = \{0\}.
    \end{multline*}
    Therefore $G^{n + 3m}$ is stable by \cref{cor:NullvarsToIntersectKernels_Flat}. The fact that $G^{n + 3m}$ is simple and irreducible follows from the fact that $G^n$ is irreducible and that $\alpha_n$ and $\beta_n$ are simple billiard loops that intersect each other only at their basepoints.
    
    To complete the proof, it remains to apply \cref{prop:Smoothing,prop:TransferGeodesicBouquet}.
\end{proof}

\bibliographystyle{alpha}
\bibliography{References}

\end{document}